\documentclass[a4paper,12pt]{article}
\usepackage{latexsym}
\usepackage{amsfonts}
\usepackage{amsmath}
\def \E {{\rm E }}
\def \R {{\rm I\!R}}

\def \N {{\rm I\!N}}
\def \t {{\tau}}
\def \v {{ v_{\epsilon, \epsilon'}}}
\def \j {{J_{\epsilon}}}
\def \k {{k(\epsilon,\epsilon')}}
\def \f {{ f_{\epsilon, \epsilon'}}}
\def \y {{ y^{\epsilon, n}}}
\def \z {{ z^{\epsilon, n}}}
\def \w {{ \dot{w}}}
\def \h {{ \dot{h}}}

\def \D {{\rm I\!D}}
\def \PrL {\left< }
\def \PrR { \right> }
\def \Quad1 {\left \! <\left<}
\def \Quad2 {\right> \! \right>}

\title{Stochastic Hyperbolic Systems, Small Perturbations and Pathwise Approximation\footnotemark[1]} 

\author{ Adnan Aboulalaa}
\date{}
\begin{document}  
\newtheorem{thm}{Theorem}[section]
\newtheorem{defi}[thm]{Definition}
\newtheorem{prop}[thm]{Proposition}
\newtheorem{Rq}[thm]{Remark}
\newtheorem{lm}{Lemma}[section]
\newtheorem{co}[thm]{Corollary}
\maketitle
\footnotetext[1]{ 
Osaka J. Math 57 (2020).  E-mail of the author : adnan.aboulalaa@polytechnique.org }

\begin{abstract}
This paper is devoted to the study of hyperbolic systems of linear partial differential equations perturbed by a Brownian motion. 
The existence and uniqueness of solutions are proved by an energy method. The specific features of this class of stochastic partial differential equations are highlighted and 
the comparison with standard existence results for SPDEs is discussed. 
The small perturbations problem is studied and a large deviation principle is stated. A pathwise approximation result, similar to the stochastic differential equations case, is established, with an application to a support theorem. 
\end{abstract}

%%%%%%%%%%%%%


%%%%%%%%%%%%%%%%%%%%%%%%%%%%%%%%%%%%%%%%%%%%%%%%%%%%%%%%%%%%%%%%%%%%%%%%
\begin{center}
\textsc{Contents}\\
\end{center}
1. Introduction\\
2. Hyperbolic stochastic partial differential systems\\
2.1. Notations\\
2.2. Remarks on the infinite dimensional stochastic calculus\\
2.2. Existence and uniqueness of a solution\\
2.3. The case of differential operators and other remarks \\
2.4. Comparison with other existence results for SPDEs\\
3. Small perturbations\\
3.1. Introduction and preliminaries\\
3.2. A large deviation principle\\
4. Pathwise approximation and applications\\
4.1 Wong-Zakai type approximation\\
4.2. Application to a support theorem\\
4.3. Application to the random evolution operator\\
5. On the regularity of the laws of the solutions\\
6. Appendix: Proofs of some technical lemmas

%%%%%%%%%%%%%%%%%%%%%%%%%%%%%%%%%%%%%%%%%%%%%%%%%%%%%%%%%%%%%%%%%%%%%%%%
\section{Introduction}
In this paper we are interested in the following class of stochastic partial differential equations:
\begin{equation}
\nonumber
(E)\left\{ \begin{array}{l}
  \displaystyle
 du(t)=\sum_{i=1}^{n}[a_{i}(t,x,D)u(t)\circ dw_{i}(t)+f_{i}(t)\circ dw_{i}(t)]+b(t,x,D)u(t)dt+g(t)dt,\\
\displaystyle
 u(0)=u_{0}\in (H^{s}(\R^{d}))^{d'},
\end{array}\right.
\end{equation}
where $H^{s}(\R^{d})$ is a Sobolev space $(s\in \R, d, d'\geq 1 )$, 
$a_{i}(t,x,D),b(t,x,D)$ are smooth families of $d'\times d'$-matrices of first
 order pseudodifferential operators (PDO), $w_{i}(t), t\in I$ are standard Wiener processes, $f_{i},g$ are continuous, possibly random, functions from 
$I=[0,T], T>0$ to $(H^{s}(\R^{d}))^{d'})$, and $\circ$ corresponds to the Fisk-Stratonovich integral or differential. Equation (E) is to be 
viewed 
as a random perturbation of a (deterministic) linear symmetric system ($a\equiv
0, f\equiv 0$, cf. Friedrichs\cite{Friedrichs}, Lax\cite{Lax}).
These systems occur often in applications, for example, the wave equation in  
non-homogeneous
media (and more generally any second order linear hyperbolic equation)
can be represented by such a system. The Maxwell equations form a symmetric hyperbolic system. 
The Dirac equation for a relativistic particle of spin 1/2 (whose
wave function is a 4-dimensional vector) is also a linear symmetric system.

Hyperbolic partial differential equations (PDE) and systems form an important class in the theory of PDE, 
and there is already a substantial literature on hyperbolic stochastic partial differential equations (SPDEs), 
see, e.g., \cite{Wa}, \cite{Hajek}, \cite{Chow}, \cite{DQS}, \cite{LR}, \cite{Kim}, \cite{Gaveau}, \cite{ACS}, \cite{LPS2}, \cite{BVW}.
These works address various models including the stochastic wave and related equations with space-time white noise, non linear models related to
conservation laws and hyperbolic systems with additive random perturbation (with assumptions on the diffusion coefficient that exclude the case
of the system ($E$)). 

Particular forms of Eq. (E) have been considered by 
Ogawa \cite{Og}, and Funaki \cite{Fu} who used a method of characteristics
to construct the solutions. Later, Kunita \cite{K2} made use of this method 
and the theory of stochastic flows to obtain solutions to nonlinear
first order partial differential equations in the scalar case. 
However, it 
is well known that, in the deterministic case, the method of characteristics 
no longer applies if $u$ is 
not scalar (and for more than one space variable). Instead, energy or 
semigroup methods are used to solve linear symmetric systems, see, e.g.,
Cordes \cite{Cordes}, H\"ormander \cite{H}, Taylor \cite{Ta}. In the first part of this paper we use 
an energy method, based on {\it a priori} estimates, to solve systems like (E). We also discuss the relation to
standard existence results for stochastic partial differential equations (SPDEs).

The second part deals with the small random perturbations of the deterministic system $ (E_{d}): du(t)=b(t,x,D)u(t)dt,\; u(0)=u_{0}$.
More precisely, if we consider the solution $u^{\epsilon} (\cdot)$ to:
\[ ({\cal S}_{\epsilon}): \; \;  du^{\epsilon}(t)=\sqrt{\epsilon}\sum_{i=1}^{n}a_{i}(t,x,D) u^{\epsilon}(t)\circ dw_{i}(t)+
b(t,x,D)u^{\epsilon}(t)dt,\; \; u^{\epsilon}(0)=u_{0},\]
\noindent
then it is easily shown that $u^{\epsilon}$ converges in probability (in some Hilbert space) to the solution $u$ of $(E_{d})$ 
as $\epsilon \rightarrow 0$. The objective is to obtain information on how the stochastic solution is close to the deterministic one as the perturbation becomes very small.
In the case of stochastic differential equations (SDEs), the deviation probabilities of the stochastic solutions from the deterministic ones, converge to $0$
 with an exponential rate (Freidlin-Wentzell estimates for SDEs), i.e., we have a large deviation principle for the small random perturbation of ordinary differential equations.
In the case of SPDEs, we face two problems: the dimension is infinite and the coefficients
 are unbounded. In the literature, these problems have been addressed
for various SPDEs, see, e.g., \cite{DZ}, \cite{Pe}, \cite{BDM}, \cite{RZ} and the references therein. In our case, it turns out that the 
Freidlin-Wentzell estimates 
hold, but with the loss of two derivatives: $u^{\epsilon}$ is
an $H^{s}$-valued process, the large deviation principle is valid
in the topology of $C^{0}(I,H^{s-2})$ associated with the norm 
$\sup_{t}|v(t)|_{s-2}$.

The third part of this paper is concerned with the approximation of the
stochastic system (E) by deterministic systems depending on a random
parameter. In the case of finite dimensional stochastic equations, this 
is sometimes called the Wong-Zakai or Stroock-Varadhan approximation \cite{WZ}, \cite{SV}. 
Such approximations have received some attention in the case of SDEs for
they allow the transfer of some properties of ordinary differential equations
to the stochastic case (e.g., the construction of the stochastic flows of
diffeomorphisms for SDEs, the initial approach to the stochastic
calculus of variations, or the approximation of the solutions), see e.g., Malliavin \cite{Ma}, Ikeda-Watanabe \cite{IW}
and the references given there. The case of SPDEs has also been considered for several models, 
see, e.g., Gy\"ongy \cite{Gy}, Brzezniak-Flandoli \cite{BF}, Hairer and Pardoux \cite{HP}, Yastrzhembskiy \cite{Yastrzhembskiy} and
Roth \cite{Roth}; the later reference considered stochastic hyperbolic equations similar to (E) in the scalar case.
In this part, the convergence of the deterministic systems associated to (E) to the stochastic one is shown and this is applied to
an infinite-dimensional extension of the Stroock-Varadhan support theorem and to 
the construction of the random evolution semigroup of Eq. ($E$).

The last part of the paper deals with the regularity of the law of the solutions to Eq.(E) using the Malliavin calculus 
techniques.

In this paper we restrict ourselves to finite dimensional Wiener processes in order to highlight the specific features of
hyperbolic systems, in particular the necessity to use Fisk-Stratonovich integral and the fact that the techniques used in parabolic SPDEs are not suitable in this case.

Finally, as we consider a class of SPDEs that involves pseudodifferential operators (PDO), let us mention some of the few works so far published which use these operators in probabilistic models. 
Among the references, we quote Kotelenez \cite{Kot}, Kallianpur and Xiong \cite{KX} and Tindel \cite{Ti} who considered
SPDEs with pseudodifferential operators with an assumption on the order ($m > d$) and a space-time white noise;
Jacob, Potrykus and Wu studied in \cite{JPW} the solution of a stochastic Burger equation driven by a space time white noise and using a PDO; and in \cite{LX}, Liu and Zhang started a study of stochastic pseudodifferential operators with a Calder\'on-type uniqueness theorem as an application to SPDEs.

%%%%%%%%%%%%%%%%%%%%%%%%%%%%%%%%%%%%%%%%%%%%%%%%%%%%%%%%%%%%%%%%%%%%%%%%

\section{Hyperbolic stochastic partial differential systems}
\subsection{Notations}
Let $d,d'\geq 1$. We denote by $S^{m}$ the set of symbols $a(x,\xi)$ of order 
$m$ on $\R^{d}$, i.e. $a\in C^{\infty}(\R^{d}\times \R^{d})$ and
for all $\alpha, \beta \in \N^{d}$ there is a constant $C(\alpha,\beta)$ 
such that $|D_{\xi}^{\alpha}D_{x}^{\beta}a(x,\xi)|\leq C(\alpha, \beta)
(1+|\xi|)^{m-|\alpha|}$, where we have used the notation: 
for $\alpha = (\alpha_{1}, \cdots, \alpha_{d}) \in \N^{d}$, $|\alpha|= \sum_{j=1}^{d}\alpha_{j}$ and $D^{\alpha}= D^{\alpha_{1}}_{1} \cdots D^{\alpha_{d}}_{d}$, where
$D^{\alpha_{j}}_{j} = (-i)^{\alpha_{j}}\partial_{j}^{\alpha_{j}}$ with $i=\sqrt{-1}$ and $\partial_{j}=\partial_{x_{j}}$, for $x=(x_{1}, \cdots, x_{d})$. 
\\
\noindent
For $a\in S^{m}$, $a(x,D)$ will denote the associated
pseudodifferential operator defined by $a(x,D)u(x)=\int a(x,\xi) \hat{u}(\xi)
e^{i x\cdot\xi}d\xi $ for $u\in C^{\infty}_{0}(\R^{d})$. ${\rm OPS^{m}}$ will
designate the set of such operators and $a^{*}(x,D)$ is the adjoint of
 $a(x,D)$. We denote by $\left< \cdot, \cdot \right>_{s}$ the scalar product on the Sobolev
space $H^{s}:=H^{s}(\R^{d}), s\in \R$. We use the same notation for the scalar
product on $(H^{s})^{d'}$.
\noindent
In the following, we shall consider matrices of pseudodifferential 
operators $a(x,D)$, which means that $a(x,D)=(a^{ij}(x,D), i,j=1,...,d')$
with $a^{ij}(x,D)\in {\rm OPS}^{m}$ for some $m$. We still denote by ${\rm OPS}^{m}$ the
set of such matrices of operators.
\\
\noindent
Throughout this paper, we fix $T>0$ and we assume that we are given a one-dimensional
 Brownian motion $w(t), t\in I:=[0,T]$ defined on a filtered probability
space $(\Omega,{\cal F}, {\cal F}_{t}, P)$ with ${\cal F}_{t}=\sigma(w(\t), 
\t\leq t)$. For $t, t' \in\R$, we set $t \wedge t'=\inf(t, t')$.  
In the sequel, we deal with $(H^{s})^{d'}$-valued process. For $p>0$, $M^{s}_{p}(I, (H^{s})^{d'})$ will designate the 
set of adapted $(H^{s})^{d'}$-valued processes $u(t),t\in I$ such that 
$\|u\|_{s,p}:=(E\sup_{t\leq T}|u(t)|_{s}^{p})^{1/p}< +\infty$. $M^{s}_{p}$ 
will generally be endowed with the norm $\|.\|_{s,p}$.
\\
\noindent
The quadratic variation of two processes $M, N$ will be denoted by $\left< M, N \right>_{t}$, the time variable subscript will be either 
$t$ or $\t$, while the inner product in $H^{s}$ will be denoted by $\left< u, w \right>_{s}$, in which case the subscript will always be the letter $s$, the order of the Sobolev space; the notation $\left< u, w \right>$, without subscript, refers to the inner product in $L^{2}((\R^{d})^{d'})$.

%%%%%%%%%%%%%%%%%%%%%%%%%%%%%%%%%%%%%%%%%%%%%%%%%%%%%
\subsection{Remarks on the infinite dimensional stochastic calculus}

As we are dealing with infinite dimensional processes, the following remarks concern the stochastic calculus concepts that will be used in this paper.
The Fisk-Stratonovich integral used in the hyperbolic SPDE (E) is defined by:

\[ \int_{0}^{T}u(t)\circ dw(t)=\int_{0}^{T}u(t)dw(t)+ \frac{1}{2}\int_{0}^{T}d\left< u,w  \PrR_{t}. \]
\noindent
where $\left< u,w \PrR_{t}$ is the quadratic cross variation of the processes $u(\cdot), w(\cdot)$. This integral is usually defined when 
the integrand $u(t)$ is a semimartingale; we also note that, in our case, $u (\cdot), w (\cdot)$ are in different spaces.

In the following $H$ is a separable Hilbert space and $(\Omega,{\cal F}, {\cal F}_{t}, P)$ is a filtered probability
space.
\\[12pt]
\noindent
$\bullet$ {\bf Martingales and semimartingales in Hilbert spaces (\cite{MP}, \cite{Metivier}, \cite{DZ}, \cite{GM}): } 
An $H$-valued process $M(t)$ is an ${\cal F}_{t}$-martingale if $\forall t, t' \geq 0$, with $t'\leq t$: (1)$M(t)$ is an ${\cal F}_{t}$-measurable random variable, 
(2)$E(\| M(t) \|) <+\infty$ and (3)$E ( M(t) | {\cal F}_{t'}) = M(t')$.
The last condition is equivalent to $E ( \left< M(t), h \PrR_{H} | {\cal F}_{t'}) = \left< M(t'), h \PrR_{H}, \forall h\in H$. 
\noindent
In this case $\| M(t) \|^{2}$ is a real submartingale and therefore $\E \| M(t) \|^{2} \leq  \E \| M(T)\|^{2}$, for all $t$ in $[0,T]$. 
We also have the following inequalities:
\begin{equation}
\label{mart1}
{\rm E} (\sup_{t\in I} \| M(t) \|^{p}) \leq (\frac{p}{p-1})^{p} \sup_{t\in I} {\rm E} \| M(t) \|^{p}, \; \; \forall p\geq 1
\end{equation}

\begin{equation}
\label{mart2}
{ \rm Pr} (\sup_{t\in I} \| M(t) \| \geq \lambda) \leq \lambda^{-p} \sup_{t\in I} {\rm E} \| M(t) \|^{p}, \; \; \forall p  > 1
\end{equation}
The space ${\cal M}^{2}(H)$ of square integrable martingales equipped with the norm:
\[ \| M(t) \|_{{\cal M}^{2}(H)}= (E (\sup_{t\in I} \| M(t) \|^{2})^{1/2)} \]
is a Banach space. $M(t)$ is said to be a local martingale if there exists an increasing sequence of stopping times $\sigma_{n}$ with $\lim \sigma_{n} = +\infty$ such that for every $n$, $M(t \wedge\sigma_{n})$ is a martingale.
\\
An $H$-valued process $X(t)$ is a semimartingale if $X(t)=M(t)+V(t)$ where $M(t)$ is a local martingale and $V(t)$ is a process which has 
a finite variation a.e. in every bounded interval. 
\\[12pt]
$\bullet$ {\bf The quadratic variation and the cross variation tensor of a martingale (\cite{MP}, \cite{Metivier}, \cite{DZ}): }
For a martingale $M(t) \in {\cal M}^{2}(H)$ there exists a unique $H\!\otimes\! H$-valued process, $ \left<\!\left< M(t) \right>\!\right> $, called the quadratic variation
of $M(t)$ such that $ M(t)\otimes M(t) - \left<\!\left<  M(t) \right>\!\right> $ is a martingale; as usual $\otimes$ denotes the tensor product of Hilbert spaces or elements
of Hilbert spaces. In the same way for two square integrable martingales $M_{1}(t), M_{2}(t)$ with 
values in two Hilbert spaces $H_{1}, H_{2}$, there exists a unique $ H_{1}\!\otimes\! H_{2}$-valued process, $ \left<\!\left< M_{1}(t), M_{2}(t)\right>\!\right> $, called the quadratic cross
variation of $M_{1}(t), M_{2}(t)$ such that $ M_{1}(t)\!\otimes\! M_{2}(t) - \left<\!\left< M_{1}(t), M_{2}(t)\right>\!\right> $ is a martingale.
\\
If $e_{i}$ and $f_{j}$ are orthonormal bases (ONB) of $H_{1}, H_{2}$, then, with $ M_{1}(t)=\sum_{i=1}^{\infty} M_{1,i}(t) e_{i}$ and 
 $M_{2}(t)=\sum_{j=1}^{\infty} M_{2,j}(t) f_{j}$, we have:
\begin{equation}
\label{Crossvariation}
\left<\!\left<  M_{1}(t), M_{2}(t)  \right>\!\right> = \sum_{i,j=1}^{\infty} \left<\!\left<  M_{1,i}(t), M_{2,j}(t)\right>\!\right>  e_{i}\!\otimes\! f_{j},  
\end{equation}
and $\left<\!\left<M_{1,i}(t), M_{2,j}(t)\right>\!\right> $ is the usual cross variation of the real martingales $M_{1,i}, M_{2,j}$. Let us note that these notions can be defined in another
way, for more general processes $X$, $Y$ with values in the Hilbert spaces $H_{1}, H_{2}$ by:
\begin{align*}
 D^{n}_{t}(X,Y) & = \sum_{k\geq 1} (X(\sigma_{n,k+1}\wedge t) - X(\sigma_{n,k} \wedge t))\!\otimes\! (Y(\sigma_{n,k+1}\wedge t) - Y(\sigma_{n,k}\wedge t)) \\
 \left<\!\left<X, Y\right>\!\right>_{t}    & = \lim_{n\longrightarrow \infty} D^{n}_{t}(X,Y), 
\end{align*}
\noindent
if the limit in probability of $ D^{n}_{t}$ exists. Here, $\sigma_{n,k}$ is family of stopping times such that $\sup_{k} (\sigma_{n,k+1}-\sigma_{n,k}) \longrightarrow 0$ as $n\longrightarrow\infty$, so that we have a subdivision of the interval $[0,t]$, see \cite{Metivier}, \cite{MP}.
\\[12pt]
\noindent
$\bullet$ {\bf It\^o  integral and It\^o formula in Hilbert spaces (\cite{DZ}, \cite{GM})}.
The It\^o integral $\int Y(t) dX(t)$ can be defined in a general setting, where $X(t)$ is a $J$-valued process, $Y(t)$ is 
a $L(J,K)$-valued process, $J,K$ are separable Banach spaces, $L(J,K)$ is the space of linear operators (not necessarily bounded) from $J$ into $K$,
with the corresponding properties, in particular the It\^o formula, see M\'etivier-Pellaumail \cite{MP}. Since we deal with the special case where $Y(t)$ belongs to a Hilbert space and $X(t)$ is a real Brownian motion, we restrict ourselves to the case of stochastic integral with respect to a $Q$-Wiener process.
\\
Let $H, U$ be two separable Hilbert spaces, $Q\in L(U)$ a linear bounded operator of $U$ which is self-adjoint, non-negative and of finite trace. A $Q$-Wiener process
is an $U$-valued process $W(t)$ such that: (1)$ W(0)=0$, (2) for all $0\leq t_{1} \leq ... \leq t_{n} \leq T$, $W(t_{1}), W(t_{2})-W(t_{1}), ..., W(t_{n})-W(t_{n-1})$ 
are independent random variables (rv's), (3)$W(t)-W(t')$ is a Gaussian random variable with the covariance operator $(t-t')Q$ and (4)the paths of $W$ are continuous $P$-a.e. If $e_{i}$ is an ONB diagonalizing $Q$ with eigenvalues $\lambda_{i}$, then $W(t)$ can be written as:
\begin{equation}
\label{QWinener}
W(t)= \sum_{i\geq 1} \lambda_{i}^{1/2} W_{i}(t) e_{i}
\end{equation}
where $W_{i}(t)$ are independent real Brownian motions. Now, let $U_{0}$ be the subspace $Q^{1/2}(U)$ of $U$; then the stochastic integral $ \int_{0}^{T} \Phi(\t) dW(\t)$ can be defined for all $L(U_{0}, H)$-valued predictable processes $\Phi(\t)$ that are Hilbert-Schmidt for all $\t$ and such that $ \int_{0}^{T} \| \Phi(\t) \|_{HS} d\t  < +\infty $.\\
\noindent
In this case, if $\phi(\t)$ is an adapted $H$-valued process and $X(t)$ is the $H$-valued process defined by:
\[X(t) = X(0) + \int_{0}^{t} \phi(\t)ds   + \int_{0}^{t} \Phi(\t) dW_{\t}, \] 
\noindent
with $X(0)$ an ${\cal F}_{0}$-measurable rv, then the It\^o formula has the following form (\cite{DZ}, \cite{GM}):
\begin{align}
\label{ItoFormula}
F(t,X(t)) =& F(0,X(0)) + \int_{0}^{t} DF(\t, F(X(\t))).\Phi(\t) dW(\t) \nonumber \\ 
          &+ \int_{0}^{t}[ \frac{\partial F}{\partial \t} (\t, X(\t))+ DF(\t, F(X(\t))).\phi(\t) \nonumber \\ 
					&+ \frac{1}{2} {\rm tr}( D^{2}F(\t,X(\t)).(\Phi(\t)Q^{1/2})(\Phi(s)Q^{1/2})^{*}]d\t, 
\end{align}
\noindent
(we use the notation $DF(x).V$ for the image of a vector $V$ by the differential map $DF(x)$). Other forms of the It\^o formula may be found in more general settings in \cite{MP}, \cite{Metivier}.

$\bullet$ {\bf Fisk-Stratonovich integral in Hilbert spaces: }
In the real case, the Fisk-Stratonovich integral with respect to a semimartingale $X(t)$ is defined by:
\begin{equation}
\label{FS}
\int_{0}^{T}Y(t)\circ dX(t)=\int_{0}^{T}Y(t)dX(t)+ \frac{1}{2}\int_{0}^{T}d\left< X,Y \PrR_{t},
\end{equation}
\noindent
where $Y(t)$ is also a semimartingale. In the last formula $X(t), Y(t)$ are also assumed to be continuous. See, e.g., \cite{IW}, \cite{Protter}.
These conditions can be relaxed; for a detailed study we refer to Meyer (\cite{Meyer}, Chap. VI, p. 109).
In the Hilbert space case, this concept seems not to have received much attention (to the best knowledge of the author); however,
the same formula (\ref{FS}) is valid with the previous definition of the cross variation (see \cite{Aida}, \cite{TN} for
some further results).
\\
In our case, for a semimartingale $u(t)=\sum_{i\geq 1} u_{i}(t) e_{i} \in H $, where the $u_{i}(t)$ are real semimartingales, and a real Brownian motion $w(t)$, we note that:
$$\left< u,w \PrR_{t}= \sum_{i\geq 1} \left< u_{i},w \PrR_{t} e_{i},$$
and by expansion in the ONB $(e_{i})$, the calculus rules in this specific Hilbert space setting can also be deduced from their real counterpart. 
For instance, if $du(t)=\sigma(u(t)) \circ dw(t)$ where $\sigma : H \longrightarrow H$ is a differentiable map, then:
\[ du(t)=\sigma(u(t)) dw(t) + \frac{1}{2} D\sigma (u(t))\cdot\sigma(u(t)) dt \]

$\bullet$ {\it Application of the It\^ o formula:} Let us consider a special case of the equation (E): $du(t)=  a u(t)\circ dw(t)$, where $u(t)$ is an adapted $(H^{s})^{d'}$-valued process, $a$ is a first order symmetric operator in $(H^{s})^{d'}$ (denoted by $H$) and $w(t)$ is a real Brownian motion. This equation is then written in the It\^ o form as:
\[ du(t)=  a u(t) dw(t) + \frac{1}{2} a^{2} u(t) dt \]
\noindent
Now we shall apply the previous It\^ o formula to $F(u(t))= |u(t)|^{2}_{s}=\left< u(t),u(t) \PrR_{s}$. For the differentials of $F$, we
have: 
$$DF (v).h= \left<  v,h \PrR_{s}+\left<  h, v \PrR_{s} \; {\rm and} \; \; D^{2}F (v).h.k= \left<  k,h \PrR_{s}+\left<  h, k \PrR_{s} $$
In our case, the $Q$-Wiener process is just a one-dimensional Brownian motion, $U=U_{0}=\R$, $Q=Id$ and $\Phi(\t)$ is the linear operator from $\R$ to $H$ defined by $\Phi(\t).\alpha = \alpha a u(\t)\in H, \alpha\in \R$ and $\Phi(\t)^{*}=\Phi(\t)$, so that: \\
\\
(a) The second order term involving $D^{2}F$ in (\ref{ItoFormula}) is:
\[ D^{2}F(u(t)).au.au=2\left<  au(t),au(t) \PrR_{s} = 2\left<  a^{*}a u(t),u(t) \PrR_{s} \]
\\
(b) For the first order derivative term $DF(\t, F(X(\t))).\phi(\t)$, we have $\phi(\t) = 1/2(a^{2}u(\t))$ and:
\begin{align*}
 DF(\t, F(X(\t))).\phi(\t) &= \frac{1}{2} ( \left<  u(\t),a^{2}u(\t) \PrR_{s}+\left< a^{2}u(\t), (\t) \PrR_{s}) \\ 
                         &= \frac{1}{2} \left<  (a^{2}+a^{*2}) u(\t), u(\t) \PrR_{s}
\end{align*}
\\
(c) The term involving $dw(\t)$ is: 
$$DF(u(\t)).a u(\t)=\left<  u(\t),a u(\t) \PrR_{\t}+\left< a u(\t), u(\t) \PrR_{s}=\left< (a+a^{*}) u(\t), u(\t) \PrR_{s}.$$
This yields:
\begin{align*}
|u(t)|^{2} =& |u(0)|^{2} + \int_{0}^{t} \left< (a+a^{*}) u(\t), u(\t) \PrR_{s} dw(\t)  \\ 
           & + \frac{1}{2}\int_{0}^{t}\left<  (a^{2}+a^{*2}) u(\t), u(\t) \PrR_{s} d\t + \frac{1}{2} \int_{0}^{t} 2\left<  a^{*}au(\t),u(\t) \PrR_{s} d\t\\
					 = & |u(0)|^{2} + \int_{0}^{t} \left< (a+a^{*}) u(\t), u(\t) \PrR_{s} dw(\t) \\
						&+ \frac{1}{2}\int_{0}^{t} \left<  (a^{*}(a^{*}+a)+ (a^{*}+a)a )u(\t),u(\t) \PrR_{s} d\t
\end{align*}
So that, with $A= a+a^{*}, L= a^{*}(a^{*}+a)+ (a^{*}+a)a$, we have:
\begin{equation}
\label{ItoFormula2}
|u(t)|^{2}=|u(0)|^{2} + \int_{0}^{t} \left<  A u(\t), u(\t) \PrR_{s} dw(\t) + \frac{1}{2}\int_{0}^{t} \left<  L u(\t),u(\t) \PrR_{s} d\t
\end{equation}
\noindent
{\bf Remark.} It is convenient to derive formulas like (\ref{ItoFormula2}) by using the symbolic stochastic and differential calculus rules, and since
this will be repeatedly used throughout the paper, let us apply these rules in this simple but typical situation:
If $X(t)=A_{X}(t)+M_{X}(t)$ and $Y(t)=A_{Y}(t)+M_{Y}(t)$ are two continuous semimartingales where $A$ and $M$ denote respectively
 the finite variation and martingale parts, then the cross variation $d\left< M_{X},M_{Y} \PrR_{t}$ is formally denoted by $dX(t)\cdot dY(t)$ and we have
a set of symbolic calculus rules including (we omit here the time parameter):
\begin{align*}
\bullet & \; Y \circ dX = Y dX + \frac{1}{2}dX\cdot dY, \\
\bullet & \; X \circ (dY\cdot dZ) = (X \circ dY)\cdot dZ=  X \circ (dY\cdot dZ),\\
\bullet & (XY)\circ dZ=X\circ(Y\circ dZ), \\
\bullet & \; d A_{X}\cdot dA_{Y}= d A_{X}\cdot d Y=0, \; \; dw\cdot dw= dt, \; \; dX\cdot dY \cdot dZ=0. 
\end{align*}
See \cite{IW}, pp.99-100. Besides these identities, the usual differential calculus rules can be applied when using the Fisk-Stratonovich differentials (as remarked before, in the present setting with Hilbert space valued processes and finite dimensional Brownian motion, these results are valid under the same conditions of finite dimension processes, via an expansion in an ONB). In our example, we apply this to $ |u(t)|_{s}^{2}= \left< u(t), u(t) \PrR_{s}$:
\begin{align}
\label{symbCalc4}
d |u(t)|_{s}^{2} =& \left< a u(t), u(t) \PrR_{s} \circ dw(t) + \left< u(t), a u(t) \PrR_{s} \circ dw(t) \nonumber \\
         =& \left< (a+a^{*}) u(t), u(t) \PrR_{s} \circ dw(t),
\end{align}
and with $A=a+a^{*}$:
\begin{equation}
\label{symbCalc5}
\left<A u(t), u(t) \PrR_{s} \circ dw(t) = \left<A u(t), u(t) \PrR_{s} dw(t)+\frac{1}{2} d \left<A u(t), u(t) \PrR_{s}\cdot dw(t).
\end{equation}
So we have to calculate cross variation $d \left<A u(t), u(t) \PrR_{s}\cdot dw(t)$; we use the above mentioned stochastic rules to get:
\begin{align}
\label{symbCalc6}
d\left<A u(t), u(t) \PrR_{s}\cdot dw(t) =& (\left<A a u(t), u(t) \PrR_{s} \circ dw(t)) \cdot dw(t)+(\left<A u(t), a u(t) \PrR_{s} \circ dw(t)) \cdot dw(t) \nonumber \\
      =& (\left<(A a  + a^{*}A) u(t), u(t) \PrR_{s} \circ dw(t)) \cdot dw(t) \nonumber \\
			 =& \left<L u(t), u(t) \PrR_{s} dw(t)) \cdot dw(t) \nonumber\\
			=& \left<L u(t), u(t) \PrR_{s} dt
\end{align}
and by replacing (\ref{symbCalc5}) and (\ref{symbCalc6}) in (\ref{symbCalc4}) we get (\ref{ItoFormula2}).

\noindent
From this calculus we may expect that the existence of solutions of systems like $(E)$ will be related to assumptions about the boundedness of
the operators $A, L$; this is possible in the case of symmetric hyperbolic systems (as in the deterministic case). We also note that this is possible only
when (E) is written in the Stratonovich form, see the remark (c) in $\S$ 2.4.
\\
%%%%%%%%%%%%%%%%%%%%%%%%%%%%%%%%%%%%%%
\subsection{Existence and uniqueness of a solution}
We consider stochastic equations of the type:
\[
({\cal E}):
\left\{ \begin{array}{l}
\displaystyle \frac{\partial u(t,x)}{\partial t}=a_{t}(x,D)u(t,x) \circ \frac{dw(t)}{dt}+b_{t}( x,D)u(t,x) +f(t,x)\circ \frac{dw(t)}{dt} +g(t,x)\\
u(0, \cdot)=u_{0} (\cdot).
\end{array}\right. \]
where $u_{0}\in (H^{s})^{d'}, f,g \in C^{0}(I, (H^{s})^{d'}) $ for some $s$, $w(t), t\in I$ is a one-dimensional Brownian motion and $a_{t}, b_{t}, t\in I$ are (matrices 
of) pseudo-differential operators satisfying the following conditions:

(i) $a_{t}(x,D), b_{t}(x,D)$ form a bounded family in ${\rm OPS}^{1}$ (see below).

(ii) The $C^{\infty}(\R^{2d})$-valued maps: $t\mapsto a_{t}^{i,j}, 
b_{t}^{i,j}$ are continuous  ($C^{\infty}(\R^{2d})$ is equipped with the usual topology of the semi-norms $\max_{x\in \R^{d}} |\partial^{p}f|, p=(p_{1}, ...,p_{k})\in \N^{k}, k\geq 0$).
\\
Throughout this paper, the space variable $x$ may not be shown and we shall write $v(t,\cdot)=v(t) \in (H^{s})^{d'}$ for $v=u, f, g$, etc.
\\[12pt]
%%%%%%%%%%  %%%%%%%%%%
{\bf Bounded family of symbols and operators:} We shall need some results about the action of pseudodifferential operators on Sobolev spaces and their boundedness: 
\\
$\bullet$ A family of operators $p_{t}(x,D)$ is said to be bounded in $H^{s}$ if: $\sup_{t\in I} |\!|\!| p_{t}(x,D) |\!|\!|_{s} \leq C$, where $|\!|\!| P|\!|\!| = \sup_{\|u\|=1} \| Pu\|$ is the
norm of the operator $P$ acting on some vector space. 
\\
$\bullet$ A family of symbols $p_{t}(x,\xi)\in S^{m}$ that satisfy $|D_{\xi}^{\alpha}D_{x}^{\beta}p_{t}(x,\xi)|\leq C_{t}(\alpha, \beta)
(1+|\xi|)^{m-|\alpha|}$ is bounded if the constants $C_{t}$ are bounded w.r.t the parameter $t$. In other words:
\begin{equation}
\label{BoundedSymbols}
\forall l > 0, \exists C_{l} >0 : \; |p_{t}|_{l}^{(m)}:= \max_{|\alpha+\beta|\leq l} \sup_{x,\xi}(D_{\xi}^{\alpha}D_{x}^{\beta}p_{t}(x,\xi) (1+|\xi|)^{-(m-|\alpha|)}) \leq C_{l} 
\end{equation}
Let us recall the important result about the action of pseudodifferntial operators on the Sobolev spaces $H^{s}$ (H\"ormander, Kumano-Go, see \cite{Kumano-Go1} (Lemma 2.1), \cite{Kumano-Go2}, p.124): for all $s>0$ there exist $C_{s}, l_{s} >0$ such that for $p(x,\xi)\in S^{m}$ we have:
\begin{equation}
\label{BoundedSymbols2}
\forall u\in H^{s} : \; \| p(x,D)u \|_{s} \leq C_{s,m} |p_{t}|_{l_{s}}^{(m)} \| u \|_{s+m}
\end{equation}
\noindent
which means that $p(x,D)$ maps $H^{s+m}$ onto $H^{s}$. The constants $C_{s,m}, l_{s}$ are independent of the operators $p$ and therefore, if $p_{t}(x,\xi)\in S^{m}$ is a bounded family in $S^{m}$ then by (\ref{BoundedSymbols}) and (\ref{BoundedSymbols2}) we will have $\| p_{t}(x,D)u \|_{s} \leq C_{s}'\| u \|_{s+m} $ for all $t\in I$. 
\\
In particular, for $m=0$ and a bounded family $p_{t}(x,\xi)\in S^{0}$, we have  $\| p_{t}(x,D)u \|_{s} \leq C_{s}'\| u \|_{s} $ for all $t\in I$, which means that the family 
of operators $p_{t}(x,D)$ is bounded in ${\rm OPS}^{0}$: 
\begin{equation}
\label{OpBoundedness}
\sup_{t\in I} |\!|\!| p_{t}(x,D) |\!|\!|_{s} \leq C
\end{equation}
%%%%%%%%%%%%%%%%

\noindent
The equation $({\cal E})$ is to be viewed as
\begin{align*}
({\cal E}_{S}): \; \; u(t)=& u_{0} +\int_{0}^{t}a_{\tau}(x,D)u(\tau)\circ dw_{\tau}+\int_{0}^{t}b_{\tau}(x,D)u(\tau) d\tau \\
                         &+\int_{0}^{t}f(\tau) \circ dw_{\tau} + \int_{0}^{t}g(\tau)d\tau,
\end{align*}
where $\circ dw(t)$ is the Fisk-Stratonovich differential. With 
 $$g_{1}(\t)=g(\t)+ \frac{1}{2}\left< f,w \PrR_{\t},$$
$({\cal E}_{S})$ can be written in the It\^{o} form as:

\begin{align*}
({\cal E}_{I}): \; \; u(t)=& u_{0} +\int_{0}^{t}a_{\tau}(x,D)u(\tau) 
dw_{\tau}+\int_{0}^{t}b_{\tau}(x,D)u(\tau) d\tau  + \int_{0}^{t}f(\tau)dw_{\tau}\\
 &+ \frac{1}{2}\int_{0}^{t}a_{\tau}(x,D)(a_{\tau}(x,D) u(\tau)+f(\t))d\tau+ \int_{0}^{t} g_{1}(\tau)d\tau.
\end{align*}

The integrals involve the processes $a_{\tau}(x,D) u(\tau)$ in the Stratonovich form $({\cal E}_{S})$ and $a_{\tau}(x,D)a_{\tau}(x,D)u(\tau)$ in the It\^{o} form; the later integrals belong to $(H^{s-2})^{d'}$ when $u(\tau)\in (H^{s})^{d'}$, the integrals above and the equations $({\cal E}_{I})$, $({\cal E}_{S})$ are to be considered 
in $(H^{s-2})^{d'}$, but it will be seen that the solution will be actually an $(H^{s})^{d'}$-valued process.
\\
\noindent
{\it Notation.} In the sequel we will set: 
\begin{eqnarray*}
A_{t}(x,D)&=& a_{t}(x,D)+a^{*}_{t}
(x,D),\\
B_{t}(x,D)&=& b_{t}(x,D)+b^{*}_{t}(x,D), 
\\ 
L_{t}(x,D)&=& A_{t}(x,D)a_{t}(x,D)+
a^{*}_{t}(x,D)A_{t}(x,D)
\end{eqnarray*}
\noindent
and we will often denote the matrices of PDOs $a_{t}(x,D), b_{t}(x,D), A_{t}(x,D),$ etc.\ by 
$a(t), b(t), A(t),$ or $a_{t}, b_{t}, A_{t},$ etc. We have the following existence and uniqueness result:
\begin{thm}
Let $a_{t}(x,D), b_{t}(x,D)$ be two families of matrices of pseudodifferential
operators which satisfy (i)-(ii). We assume that $\sup_{t}
E|f(t)|_{s+1}^{4}, \sup_{t}E|g_{1}(t)|_{s}^{4}$ are bounded, with $g_{1}(t)=g(t)+ \left< f,w \PrR_{t}$, and\\
{\rm (iii)} $A_{t}(x,D), B_{t}(x,D)$ are bounded families in ${\rm OPS}^{0}$\\
{\rm (iv)} $L_{t}(x,D)$ is a bounded family in ${\rm OPS}^{0}$\\
Then $({\cal E})$ has a unique solution $u\in M_{2}(I, H^{s})$ and 
$(u(t)), t\in I$ is a strong Markov process. 
Moreover, the solution $u(t)$, viewed as an $H^{s'}$-valued process ($s'\leq s$), has a modification which is almost-surely $\gamma$-H\"older 
continuous with respect to the norm $|\cdot|_{s-2}$, for all $\gamma \in ]0,1/4[$.
\end{thm}

This theorem extends easily to the case of symmetrizable systems: instead of
conditions (iii) and (iv), we may suppose that there exist smooth families of 
$(d'\times d')$-matrices of pseudodifferential operators $R_{t}^{1},
R_{t}^{2} \in {\rm OPS}^{0}$
such that the principal symbols $R_{0}^{i}(t,x,\xi), i=1,2$ are positive 
definite 
matrices for $|\xi|\geq 1$, and (iii), (iv) are to be replaced by:\\
(iii)  $R_{t}^{1}a_{t}+a^{*}_{t}R_{t}^{1}, R_{t}^{2}b_{t}+b^{*}_{t}R_{t}^{2} $
 form bounded families in ${\rm OPS}^{0}$.\\
\noindent
(iv)  $(R_{t}^{1}a_{t}+a^{*}_{t}R_{t}^{1})a_{t}+a^{*}_{t}
      (R_{t}^{1}a_{t}+a^{*}_{t}R_{t}^{1})$ is a bounded family in ${\rm OPS}^{0}$.

\noindent
The rest of this section is devoted to the proof of Theorem 2.1, which is inspired from the corresponding proof
for deterministic symmetric systems. We mention that similar {\it a priori} estimates were used by Mohammed and Sango \cite{MS2} to obtain existence results for an hyperbolic PDE with additive random term defined by a finite dimensional Brownian motion, and in \cite{GT}, Grecksch and Tudor used 
a parabolic regularization method to solve first order stochastic equations similar to (E) in the scalar case.
\\
\noindent
{\it 2.3.1 Energy estimates.} 
Let $u\in M_{2}(I, (H^{s})^{d'}), f, g\in C^{0}(I,(H^{s})^{d'})$ be such that
\begin{align*}
u(t)=& u_{0}+\int_{0}^{t}a_{\tau}(x,D)u(\tau)\circ dw_{\tau}+\int_{0}^{t}b_{ \tau}(x,D)u(\tau) d\tau \\
 &+ \int_{0}^{t}(f(\tau)\circ dw_{\tau} +g(\tau)d\tau).
\end{align*}
Then by the It\^o formula ($\S 2.2$) we get:
\begin{align}
\label{ItoFormula4}
|u(t)|^{2}_{s}=& |u_{0}|^{2}_{s}+\int_{0}^{t}(\left<  A(\t)u(\t),u(\t) \PrR_{s}+ 2Re \left<  f(\t), u(\t) \PrR_{s})dw(\t) \nonumber \\
&+\int_{0}^{t}\left< B(\t)u(\t),u(\t) \PrR_{s}d\t +\frac{1}{2}\int_{0}^{t} \left< L(\t)u(\t),u(\t) \PrR_{s} d\t \nonumber \\
&+Re \int_{0}^{t}\{\left< A(\t)u(\tau), f(\tau) \PrR_{s}+2\left< u(\t),g_{1}(\t) \PrR_{s} \nonumber  \\
&+\left< a(\t)u(\t), f(\t) \PrR_{s}+\left< f(\t), f(\t) \PrR_{s})d\t  \}
\end{align}
This formula can be obtained in the same way as for the case treated in $\S 2.2$ formula (\ref{ItoFormula2}) where $b, f, g$ are $\equiv 0$, which explains the main terms $A,L$.

\noindent
Using the boundedness of the family $a_{t},\; b_{t}$ and Schwarz's and martingale inequalities, it 
follows from (\ref{ItoFormula4}) that
\begin{align*}
E \sup_{\theta\leq t}|u(\theta)|_{s}^{4}\leq& C_{2}(E|u_{0}|^{4}_{s}+ E\{
\int_{0}^{t}(A^{2}|u(\t)|^{4}_{s}
+|f(\t)|_{s}^{2}|u(\t)|^{2}_{s})d\t \\
&+\int_{0}^{t}((B^{2}+L^{2})|u(\t)|^{4}_{s}+ |g_{1}(\t)|^{2}_{s}|u(\t)|^{2}_{s})
d\t\\
&+\int_{0}^{t} (A^{2}|f(\t)|^{2}_{s}|u(\t)|^{2}_{s}+a^{2}|f(\t)|_{s+1}^{2}
|u(\t)|^{2}_{s} +|f(\t)|^{4}_{s})d\t\}
\end{align*}
If we set $\phi^{2}(t)= E\sup_{\theta\leq t}|u(\theta)|^{4}$, the last 
inequality yields
\begin{align*}
\phi^{2}(t)\leq& C_{3}\{E|u_{0}|^{4}+\int_{0}^{t}[\phi^{2}(\t)+\phi(\t)(
             (E|f(\t)|_{s+1}^{4})^{1/2} \\
             &+ (E|g(\t)|_{s}^{4})^{1/2}) +E|f(\t)|_{s}^{4} ]d\t\}\\
    \leq& C\{ E|u_{0}|^{4}+\int_{0}^{T}(E|f(\t)|_{s+1}^{4}+E|g_{1}(\t)|_{s}^{4})
         d\t +\int_{0}^{t}\phi^{2}(\t)d\t \}
\end{align*}  
which implies by the Gronwall lemma the following estimate:
 \begin{equation}
\label{energy1}
E \sup_{\theta\leq T}|u(\theta)|_{s}^{4}\leq C (E|u_{0}|_{s}^{4}+\int_{0}^{T}
 (E|f(\t)|^{4}_{s+1}+E|g_{1}(\t)|^{4}_{s})d\t ).
\end{equation}
\noindent
%%%%%%%% %%%%%%%%
\begin{Rq}\label{Rk-p-estimate} {\rm We can show that $E \sup_{\theta\leq T}|u(\theta)|_{s}^{p}$ is bounded for $p\geq 1$ in a similar way:
using the moment inequality for martingales:
\begin{equation}
\label{MartingaleMoment}
E \sup_{\theta\leq t} |M_{\theta}|^{2p} \leq K_{p} E \left<\left<M\PrR\PrR_{t}^{p},
\end{equation}
we have for all $1 \leq p<\infty$ and for the following term, obtained when estimating $E \sup_{\theta\leq t} |u(\t)|_{s}^{2p}$ via 
the formula (\ref{ItoFormula4}) which expands $|u(\t)|_{s}^{2}$:
\begin{align*}
E \sup_{\theta\leq t}(\int_{0}^{\theta}\left<  A(\t)u(\t),u(\t) \PrR_{s} dw(\tau))^{2p} &\leq K_{p} {\rm E} (\int_{0}^{t}\left<  A(\t)u(\t),u(\t) \PrR_{s}^{2})^{p} d\tau \\
                                               & \leq  K_{p} T^{p-1} {\rm E} \int_{0}^{t} C_{A}^{2}\left< u(\t),u(\t) \PrR_{s}^{2p} d\tau \\
                                               & \leq K_{p} T^{p-1} C_{A}^{2} \int_{0}^{t} {\rm E} \sup_{\theta\leq \t} \left< u(\theta),u(\theta) \PrR_{s}^{2p} d\tau,
\end{align*}
where we have used H\"older's inequality. The other terms found by expanding $ {\rm E} \sup_{\theta\leq T}|u(\theta)|_{s}^{p}$ are 
treated in the same way as for $p=2$, so that by setting  $\phi^{2}(t)= {\rm E} \sup_{\theta\leq t}|u(\theta)|^{2p}$ we get an estimate of the type:
\begin{equation}
\label{energie11}
 \phi^{2}(t) \leq  {\rm E} |u_{0}|^{2p} + C_{1} + C_{2 } \int_{0}^{t} \phi^{2}(\t) d\t + C_{3} \int_{0}^{t} \phi(\t) d\t 
\end{equation}
\noindent
which implies the boundedness of ${\rm E} \sup_{\theta\leq T}|u(\theta)|^{2p}$ by the Gronwall lemma (In the following, by the Gronwall lemma we mean also its extension, cf. Bihari \cite{Bihari}; however, the use of this extension will not be necessary in general. In the previous case for instance, we can just note that $\int_{0}^{t} \phi(\t)d\t \leq 
 1+T \int_{0}^{t} \phi^{2}(\t)d\t$ and use the standard Gronwall lemma.)
}
\end{Rq}
%%%%%%%% %%%%%%%%
\noindent
{\it 2.3.2 Construction of the solution.}\\
\noindent
{\it (a) Preliminaries.} Let $\chi\in C_{0}^{\infty}(\R^{d}) $ be a 
test function with $\chi\geq 0, \chi(-x)=\chi(x)$ and $\int \chi(x)dx=1$ .
 Given $\epsilon>0$, let
$J_{\epsilon}$ be the Friedrichs mollifier defined by
\[ J_{\epsilon}(v)(x)=\int\chi_{\epsilon}(x-y)v(y)dy \; \; \mbox{for}\;
v\in L^{2}\; \; \mbox{with}\;\chi_{\epsilon}(x)=\frac{1}{\epsilon^{n}}
\chi(x/\epsilon).\]
%%%%%%% %%%%%%
\noindent
We recall the following properties:\\
$\bullet$ $J_{\epsilon}$ maps continuously $H^{s}$ into $H^{\infty}:= \bigcap_{s\geq 1} H^{s}$, equipped with the projective topology, in particular we have (see \cite{Treves}, Proposition 4.1, p.114):
\begin{equation} 
\label{Mollifier1}
 \forall \epsilon > 0  \; \exists C_{\epsilon, k} > 0 : \;   |J_{\epsilon} v|_{s+k} \leq C_{\epsilon, k} |v|_{s}
\end{equation}
$\bullet$ {\it Friedrichs lemma on commutators:} The commutators $[J_{\epsilon}, p_{t}(x,D)]$ remain in a bounded set of pseudodifferential operators of order $m-1$
if the $p_{t}(x,D), t\in I$ belong to a family of bounded pseudodifferential operators of order $m$; see \cite{Treves} Theorem 4.1. p. 116 and Remark 4.1. p. 118, see also \cite{Cordes} p. 79 and p.204 for similar results.\\
%%%%%%% %%%%%%%
The family of operators $a_{t}(x,D)J_{\epsilon}, b_{t}(x,D)J_{\epsilon}$ is bounded as a family of  $L((H^{s})^{d'})$: indeed as by the assumption (i), $a_{t}(x,D)$ is 
a bounded family in ${\rm OPS}^{1}$, we have by (\ref{BoundedSymbols2}) and (\ref{Mollifier1}):
\[ |a_{t}(x,D)J_{\epsilon} u |_{s} \leq C_{a} |J_{\epsilon} u |_{s+1} \leq   C_{a}   C_{\epsilon} | u|_{s}, \; \forall t\in I \]
\[ | J_{\epsilon}a_{t}(x,D) u |_{s} \leq C_{\epsilon} |a_{t}(x,D) u |_{s-1} \leq   C_{a}   C_{\epsilon} | u |_{s}, \; \forall t\in I \]
We have the same bounds for $ b_{t}(x,D) J_{\epsilon}, J_{\epsilon}b_{t}(x,D)$. Hence $a_{t}(x,D) J_{\epsilon}$ and $b_{t}(x,D) J_{\epsilon}$ are continuous (bounded) 
in $L((H^{s})^{d'})$ and for $\epsilon$ fixed their norm is uniformly bounded in $t$. Now let us consider the equation: 
\begin{align*}
({\cal E}_{\epsilon}): \; \; u(t)=& u_{0} +\int_{0}^{t}a_{\tau}(x,D)J_{\epsilon}u(\tau)\circ dw_{\tau}+\int_{0}^{t}b_{\tau}(x,D)J_{\epsilon}u(\tau) d\tau \\
 &+ \int_{0}^{t}(f(\tau)\circ dw_{\tau} +g(\tau)d\tau).
\end{align*}
The equation $({\cal E}_{\epsilon})$  can be viewed as an SDE in a Hilbert space: $du(t) = a_{\epsilon}(t)u(t)dw(t) + 
b_{\epsilon}(t)u(t)dt + f(t)dw(t) + g(t)dt$  with $w(t)$ a Brownian motion and
as $a_{\epsilon}(t), b_{\epsilon}(t)$ are bounded operators; the local Lipshitz property holds for the coefficients, and $({\cal E}_{\epsilon})$ has a unique solution (see e.g. M\'etivier-Pellaumail 
\cite{MP} $\S$ 6.10 p.74, M\'etivier \cite{Metivier}).
\begin{Rq} {\rm The solution to  (${\cal E}$) will be constructed as the limit of the solutions $u_{\epsilon}$ of the SDEs (${\cal E}_{\epsilon}$) in the Hilbert space $H^{s}$.
The solution $u_{\epsilon}(t)$ verifies (${\cal E}_{\epsilon}$) in the It\^o and Stratonovich forms and it has a modification which is a.e. H\"older continuous w.r.t the norm of $H^{s}$; this property is proved by the same argument that will be used for $u$ (see the \S (c.4) in the proof of Theorem 2.1).
}
\end{Rq}
%%%%%% %%%%%%%%
\noindent
{\it Notation.} We set: \\ 
$A_{\epsilon}(t)=a_{t}(x,D)J_{\epsilon}+J_{\epsilon}
a_{t}^{*}(x,D)$ and $L_{\epsilon}(t)=A(t)a_{t}(x,D)J_{\epsilon}+J_{\epsilon}
a^{*}_{t}(x,D)A(t)$.
%%%%%%%%%   %%%%%%%%
\begin{lm}
Under the conditions of Theorem 2.1, the operators $A_{\epsilon}(t), L_{\epsilon}(t)$ form a bounded family of operators in ${\rm OPS}^{0}$, that is:
for all $s >0$ these operators are continuous in $(H^{s})^{d'}$ and their norms are uniformly bounded in $t$ and $\epsilon$ : $|\!|\!| A_{\epsilon}(t) |\!|\!|_{s} \leq C_{s, A}$ and
 $|\!|\!| L_{\epsilon}(t) |\!|\!|_{s} \leq C_{s,L}$ for some constants $C_{s,A}, C_{s,L}$.
\end{lm}
{\it Proof.} For $A_{\epsilon}$ we write:
\[A_{\epsilon}(t)= J_{\epsilon} (a_{t}+a^{*}_{t})+[a_{t}, J_{\epsilon}] = A(t)J_{\epsilon}+[a_{t}, J_{\epsilon}]  \]
\noindent
By assumption (conditions (iii), (iv) of Theorem 2.1), $A(t), L(t)$ are bounded families of operators: 
$|\!|\!| A(t) |\!|\!|_{s} \leq C_{s,A}, |\!|\!| L(t) |\!|\!|_{s} \leq C_{s,L}$ and the fact that $ a^{*}_{t}$ is a bounded family in ${\rm OPS}^{1}$ implies by the Friedrichs lemma that $[a_{t}, J_{\epsilon}]$ is a bounded family in ${\rm OPS}^{0}$; 
this implies that it is bounded in $L((H^{s})^{d'})$ for all $s >0$ by the above mentioned results about the action of pseudodiffeerential operators on Sobolev spaces, see the bounds  (\ref{BoundedSymbols2}) and (\ref{OpBoundedness}). As for the family $L_{\epsilon}(t)$ we write: 
\[ L_{\epsilon}(t)= J_{\epsilon}(A(t)a_{t}+a_{t}^{*}A(t)) +[A(t)a_{t}, J_{\epsilon}] = L(t)J_{\epsilon}+[A(t)a_{t}, J_{\epsilon}] \]
\noindent
and we use the same argument: $L(t)$ is a bounded family in ${\rm OPS}^{0}$ by the assumption (iv) and $A(t)a_{t}$ is a bounded family in ${\rm OPS}^{1}$ ($A(t), a(t)$ being bounded families in ${\rm OPS}^{0}$ and ${\rm OPS}^{1}$ respectively). $\Box$

\noindent
Using this lemma we can prove as in \S 2.3.1 (estimate (\ref{energy1})) the following estimates:
\begin{equation}
\label{energy2}
E\sup_{\theta\leq t}|u_{\epsilon}|^{4}_{s} \leq  C
(E|u_{0}|_{s}^{4}+\int_{0}^{T}[E|f(\t)|^{4}_{s+1}+ E|g(\t)|^{4}_{s}]d\t),
\end{equation} 
the constant being independent of $\epsilon$.

\noindent
{\it (b) The construction.} Let $u_{0}\in (H^{s+2})^{d'},f,g\in C^{0}(I, 
(H^{s+2})^{d'})$ 
and  $u_{\epsilon}$
be the solution to $({\cal E}_{\epsilon})$ with the above data.
Let $v_{\epsilon,\epsilon'}=u_{\epsilon}-u_{\epsilon'}$. Then
\begin{align*}
d \v=& a(t,x,D)\j \v(t)\circ dw(t)+b(t,x,D)\j\v(t)dt\\
     &+ \f(t)\circ dw(t) +g_{\epsilon,\epsilon'}(t)dt,
\end{align*}
with
\[ \f(t)=a(t,x,D)(\j-J_{\epsilon'})u_{\epsilon}(t),\; g_{\epsilon,\epsilon'}(t)=b(t,x,D)(\j-J_{\epsilon'}) u_{\epsilon}(t).\]
%%%%%%% %%%%%
In order to simplify the proof, we assume that $b=0$ and then $g_{\epsilon, \epsilon'}=0$; in the proof below, these purely deterministic terms will give rise to terms 
that are similar to those which appear in the deterministic hyperbolic systems and have no interaction with the stochastic terms.
%%%%%  %%%%%
\begin{lm}
There exists $\k >0$ with $\k\rightarrow 0$ as $\epsilon,\epsilon'\rightarrow
 0$ such that for all $v\in H^{s+1}$.
\begin{equation}
 |(J_{\epsilon}-J_{\epsilon'})v|_{s}\leq \k |v|_{s+1}
\end{equation}
\end{lm}
\noindent
{\it Proof.} Observe that $\widehat{((J_{\epsilon}-J_{\epsilon'})v)}(\xi)=
(\hat{\chi}
(\epsilon \xi)-\hat{\chi}(\epsilon' \xi))\hat{v}(\xi)$, which implies that
\[  |(J_{\epsilon}-J_{\epsilon'})v|_{s}\leq \sup_{\xi\in \R^{d}}\frac{|
\hat{\chi}(\epsilon \xi)-\hat{\chi}(\epsilon' \xi)|}{(1+|\xi|^{2})^{1/2}}
 |v|_{s+1}.\]
Then the lemma follows from the fact that $\k:=\sup_{\xi\in \R^{d}}|
\hat{\chi}(\epsilon \xi)-\hat{\chi}(\epsilon' \xi)|/(1+|\xi|^{2})^{1/2}
\rightarrow 0$ as $\epsilon, \epsilon'\rightarrow 0$. $\Box$

\begin{lm}
\label{lm-Cauchy}
 $(u_{\epsilon})$ is a Cauchy family in $M_{4}(I,(H^{s})^{d'})$, namely:
\[ E\sup_{t\in I}|u_{\epsilon}(t)-u_{\epsilon'}(t)|^{4}_{s} \rightarrow 0 \; 
as \;\epsilon,\epsilon' \rightarrow 0. \]
\end{lm}
\noindent
{\it Proof}. We have 
\[ d|\v(t)|_{s}^{2}= \left<  A_{\epsilon}(t)\v(t),\v(t) \PrR_{s}\circ dw(t)
+2 Re \left<  \f(t), \v(t) \PrR_{s}\circ dw(t) \]
We write $|\v(t)|_{s}^{2}$ in the It\^o form, and using martingale and Schwarz's inequalities we get:
\begin{align*}
E\sup_{\theta\leq t} |\v(\theta)|_{s}^{4}&\leq C E\int_{0}^{t}\{ A|\v(\t)|^{4}
   + |\f(\t)|^{2}|\v(\t)|^{2}\\
&+L|\v(\t)|^{4} + |\v(\t)|^{2}(|\f(\t)|^{2}_{s+1}+|\f(\t)|^{2}_{s+2})\\
&+|\v(\t)|^{2}(|(J_{\epsilon}-J_{\epsilon'})(a_{\t}(x,D)J_{\epsilon}
    u_{\epsilon}(\t)+f(\t))|_{s+1}^{2}\\
&+|\f(\t)|^{4}_{s}\} d\t.
\end{align*}
We have $|f_{\epsilon,\epsilon'}(t)|_{s}\leq C|(\j-J_{\epsilon'})
u_{\epsilon}(t)|_{s+1}\leq \k |u_{\epsilon}(t)|_{s+2}$ by Lemma 2.2, therefore
\begin{align*}
E\sup_{\theta\leq t} |\v(\theta)|_{s}^{4}\leq& CE\int_{0}^{t}\{
|\v(\t)|_{s}^{4} + \k^{2} |v(\t)|^{2} |u_{\epsilon}(\t)|^{2}_{s+4}\\
   &+\k^{2} |\v(\t)|^{2}|f(\t)|^{2}_{s+2} + \k^{4}|u_{\epsilon}(\t)|^{4}_{s+1}. 
\end{align*}
Hence, by setting $\phi^{2}_{\epsilon,\epsilon'}(t)= E\sup_{\theta\leq t}|\v(\theta)|^{4}$ and
using the boundedness of $E\sup_{t\leq T}|u_{\epsilon}(t)|_{s}^{2}$  we get
\[\phi^{2}_{\epsilon,\epsilon'}(t) \leq C\int_{0}^{t}(
\phi^{2}_{\epsilon,\epsilon'}(\t)+\k \phi_{\epsilon, \epsilon'}(\t) +\k)d\t. \]
from which we deduce (by the Gronwall lemma) that $E\sup_{t\in I}|u_{\epsilon}(t)-u_{\epsilon'}(t)
|_{s}^{4}\rightarrow 0$ as $\epsilon,\epsilon'\rightarrow 0$. $\Box$

\begin{lm} \label{regul-limit}
Let $u$ be the limit of $(u_{\epsilon})$ as $\epsilon\rightarrow 0$. Then 
$u$ satisfies the equation (${\cal E}$).
\end{lm}
{\it Proof.} Let
\begin{align*}
F(t) =& u(t)-u(0)-\int_{0}^{t}a_{\t}(x,D) u(\t)dw(\t) -\int_{0}^{t}f(\t)dw(\t)\\
   &-\frac{1}{2} a_{\t}(x,D)a_{\t}(x,D) u(\t)d\t -\int_{0}^{t}g_{1}(\t)d\t .
\end{align*}
We want to show that $F(t)=0 $ a.e. Since $u_{\epsilon}$ is a solution to $({\cal E}_{\epsilon})$, we have:
\begin{align*}
F(t) =& u(t)-u_{\epsilon}(t)-\int_{0}^{t}a_{\t}(x,D) (u(\t)-J_{\epsilon} u(\t))
dw(\t) \\
   &-\frac{1}{2} \int_{0}^{t} a_{\t}(x,D) (a_{\t}(x,D) u(\t)- J_{\epsilon}a_{\t}(x,D) J_{\epsilon} u_{\epsilon}(\t))d\t.
\end{align*}
it follows form the boundedness of the family $(a_{t}(x,D))$ in ${\rm OPS}^{1}$ and martingale inequalities that
\begin{align*}
E|F(t)|_{s-2}^{2} \leq& C(E|u(t)-u_{\epsilon}(t)|_{s-2}^{2}+\int_{0}^{t}E|u(\t)-J_{\epsilon} u(\t)|_{s-1}^{2} d\t\\
&+\int_{0}^{t} E|u(\t)-u_{\epsilon}(\t)|_{s}^{2}d \t.
\end{align*}
Then we get $E |F(t)|_{s-2}^{2}=0$ by letting $\epsilon\rightarrow 0$ in the 
last inequality. $\Box$
\\[12pt]
\noindent
{\it (c) End of the proof of theorem 2.1:}

\noindent 
{\it (c.1)} Let $u_{0}\in (H^{s})^{d'}, f\in C^{0}(I, (H^{s+1})^{d'})$ and 
$u_{0}^{\epsilon}\in (H^{s+2})^{d'},f^{\epsilon}\in C^{0}(I, (H^{s+2})^{d'})$ 
be such that
\[ |u_{0}^{\epsilon}-u_{0}|_{s} \; \; \mbox{and}\;\; E  \sup_{t\in I}
|f^{\epsilon}(t)-f(t)|_{s+1}^{4}\longrightarrow 0, \]
as $\epsilon \rightarrow 0$. Let $u^{\epsilon}$ be the solution of $({\cal E}_{\epsilon})$
with the data $u_{0}^{\epsilon}, f^{\epsilon}$. Then $u^{\epsilon}-u^{\epsilon'}$ satisfies:
\begin{equation}
\label{end1}
 d(u^{\epsilon}(t)-u^{\epsilon'}(t))= a_{t}(x,D)(u^{\epsilon}(t)-
u^{\epsilon'}(t))\circ dw(t) + (f^{\epsilon}(t)-f^{\epsilon'}(t))\circ dw(t).
\end{equation}
Let $(\phi_{\epsilon,\epsilon'}(t))^{2}= E\sup_{\t\leq t}|u^{\epsilon}(\t)-
u(\epsilon')(\t)|^{4}_{s}$. The same calculations that give (\ref{energie11}) applied this time to 
(\ref{end1}) yields
\begin{align}
\label{end2}
(\phi_{\epsilon,\epsilon'}(t))^{2} \leq & C\int_{0}^{t}\{(\phi_{\epsilon, \epsilon'}(\t))^{2} + \phi_{\epsilon,\epsilon'}(\t) (E|f^{\epsilon}(\t)
-f^{\epsilon'}(\t)|_{s+1}^{4})^{1/2}\} d\t \nonumber \\
&+\int_{0}^{T} E|f^{\epsilon}(\t)
-f^{\epsilon'}(\t)|_{s}^{4} d\t +|u_{0}^{\epsilon}-u_{0}^{\epsilon'}|_{s}^{4}.
\end{align}
But the same energy inequality (\ref{energy1}) applied to the equation satisfied by
$u^{\epsilon}$ implies that $E\sup_{t\in I}|u^{\epsilon}(t)|^{4}_{s}$ is 
bounded by a constant independent of $\epsilon$. Therefore 
$\phi_{\epsilon,\epsilon'}(T)$ is bounded. From this and (\ref{end2}), we deduce that
\begin{align*}
(\phi_{\epsilon,\epsilon'}(t))^{2} &\leq  C\int_{0}^{t}\{(\phi_{\epsilon,\epsilon'}(\t))^{2} +\int_{0}^{T}E|f^{\epsilon}(\t)
-f^{\epsilon'}(\t)|_{s+1}^{4})^{1/2}d\t\\
&+\int_{0}^{T}E|f^{\epsilon}(\t)-f^{\epsilon'}(\t)|_{s}^{4} d\t +
|u_{0}^{\epsilon}-u_{0}^{\epsilon'}|_{s}^{4},
\end{align*}
which implies that $\phi_{\epsilon,\epsilon'}(T)\rightarrow 0$ as $\epsilon,\epsilon'\rightarrow 0$. Hence $u^{\epsilon}$ is a Cauchy family in
$M_{2}(I, H^{s})$. Finally, the fact that its limit $u$ satisfies the 
equation $du(t)=a(t,x,D)\circ dw(t)$ can be proved exactly as in Lemma 2.4.
\\
\noindent
This proves the existence of a solution to equation (${\cal E}$).
\\
\noindent
{\it (c.2)} The uniqueness follows from the energy estimate (\ref{energy1}).
\\
\noindent
{\it (c.3)} The Markov property of the process $u$ can be proved as usual 
(using the fact that $w$ is of independent increments) see, e.g., 
Da Prato-Zabczyk \cite{DZ}.
\\
\noindent
{\it (c.4)} To prove the continuity of the solution, we recall the Kolmogorov-Centsov theorem for metric-space valued processes, see Kallenberg \cite{Kallenberg}:
let $(E,d)$ be a complete metric space, and $X(t), t\in I \subset \R^{d}$ an $E$-valued process such that there exist $C, \alpha, \beta > 0$:
\begin{equation}
\label{ContinuityCriterion}
\E (d(X(t), X(t'))^{\alpha} \leq C |t - t'|^{d+\beta} \; \; \forall t, t' \in I,
\end{equation}
then $X(t)$ has a modification which is almost-surely $\gamma$-H\"older continuous for all $\gamma \in ]0, \beta/\alpha[$ and this modification verifies (\ref{ContinuityCriterion}).
To simplify the proof we suppose that $b=0$ and $f=g=0$; in this case, we have, for the solution $u$ to Eq. (${\cal E}$):
\[ u(t) - u (t') = \int_{t'}^{t}a(\t) u(\t) dw(\t) + \frac{1}{2}\int_{t'}^{t}a(\t)^{2} u(\t) d \t \]
and we apply the previous criterion for $\alpha=4$: the quantities $\E |u(t) - u (t')|_{s-2}^{4}$ are controlled by the sum of terms like:
\[  T_{ij}= \E (|\int_{t'}^{t}a(\t) u(\t) dw(\t)|_{s-2}^{j}|\int_{t'}^{t}a(\t)^{2} u(\t) d\t)|_{s-2}^{j}), \; i,j =0,...,4. \]
We use the Schwarz and moment martingale inequalities (\ref{MartingaleMoment}) to estimate these terms, for instance:
\begin{align*}
\E|\int_{t'}^{t}a(\t) u(\t) dw(\t)|_{s-2}^{4} & \leq K_{2} \E (\int_{t'}^{t} |a(\t) u(\t)|_{s-2}^{2} d\t)^{2} \\
                             & \leq K_{2} \E \sup_{\theta \in I}|a(\theta) u(\theta)|_{s-2}^{4} (\int_{t'}^{t} d\t)^{2}   \\
													&  \leq K_{2} C_{2} (t-t')^{2},
\end{align*}
where $C_{2} = \E  \sup_{\theta \in I}|a(\theta) u(\theta)|_{s-2}^{4} \leq A  \E \sup_{\theta \in I}|u(\theta)|_{s-1}^{4}$, which is bounded (Remark \ref{Rk-p-estimate}). 
For the other terms we can show in the same way that $T_{ij} \leq K_{ij} |t-t'|^{k}$, with $k\geq 2$. The estimates of the terms $T_{ij}$ which include the quantities $|a(\t)^{2} u(\t)|_{s-2}$ will be done by $\E \sup_{\theta \in I}|u(\theta)|_{s}^{k}$, this explains the 
continuity w.r.t the norm $|.|_{s-2}$. (\ref{ContinuityCriterion}) is then verified with $\alpha=4, \beta=1, d=1$ and then $\gamma \in ]0, 1/4 [$. $\Box$

%%%%%%%%%%%%%%%%%%%%%%%%%%%%%%%%%%%%%%%%%%%%%%%%%%%%%%%%%%%%%%%%%%%%%%%%%%%%%%%%%%%%%%%%%%
\subsection{The case of differential operators and other remarks}
{\bf (a) The case of hyperbolic differential systems:} In this section we consider the equation $({\cal E})$ where $a$ and 
$b$ are first order differential operators:
\[ a_{t}(x,D):= \sum_{i=1}^{d}\alpha^{i}(t,x)\frac{\partial}{\partial x^{i}} ,
  \; b_{t}(x,D):= \sum_{i=1}^{d}\beta^{i}(t,x)\frac{\partial}{\partial x^{i}} ,
\]
where $\alpha^{i}(t,x),\beta^{i}(t,x)$ are symmetric $d\times d'$-matrices; $A(t)$ and the adjoint of $a_{t}$ are given by:
\[
a_{t}^{*}u = -\sum_{i=1}^{d}\alpha^{i}(t,x)\frac{\partial u}{\partial x^{i}} - \sum_{i=1}^{d}\frac{\partial \alpha^{i}}{\partial x^{i}}u, \; 
A(t)u=- \sum_{i=1}^{d}\frac{\partial \alpha^{i}}{\partial x^{i}}u.
\]
The same formulas hold for $b$, so that the condition (iii) of Theorem 2.1 is satisfied if the first order partial derivatives of $\alpha^{i}(t,x)$ and $\beta^{i}(t,x)$ are bounded. As for $L(t)=A_{t}a_{t}+a^{*}_{t}A_{t}$, we have:
\[ 
L(t)u = \sum_{i, k=1}^{d}[- \alpha^{i}\frac{\partial^{2} \alpha^{k}}{\partial x^{i}\partial x^{k}} + \frac{\partial 
  \alpha^{k}}{\partial x^{k}}\frac{\partial \alpha^{i}}{\partial x^{i}}] u.
\]
The condition (iv) is then satisfied if the $\alpha^{i}$ and their 2d order derivatives are bounded. 
\\[12pt]
\noindent
{\bf (b) Scalar equations:} In the rest of this paragraph we focus on the scalar case ($d'=1, d=1$) for simplicity, and we consider the equation (${\cal E}_{S}$):
\begin{equation}
\nonumber
({\cal E}_{S})\left\{ \begin{array}{l}
  \displaystyle
 du(t)= a(t,x,D)u(t,x) \circ dw(t) + b(t,x,D)u(t,x) dt + f(t,x)\circ dw(t)+g(t,x) dt,\\
\displaystyle
 u(0,x)=u_{0}(x),
\end{array}\right.
\end{equation}
with
\[ a(t,x,D)u= a^{1}(t,x){\partial u}/{\partial x} + a^{0}(t,x)u, \; \;  b(t,x,D)u= b^{1}(t,x){\partial u}/{
 \partial x} + b^{0}(t,x)u. \]
\noindent
Thus, according to Theorem 2.1, for each $u_{0}\in H^{s}, f\in C^{0}(I,H^{s+1}), g\in C^{0}(I, H^{s})$, 
Eq. (${\cal E}_{S}$) has a unique solution.
\noindent
In the case of regular data $u_{0}, f, g$, Ogawa \cite{Og} and Funaki \cite{Fu} gave an
expression of the solutions to (${\cal E}$) in particular cases using 
a stochastic version of the classical characteristic method. Kunita \cite{K1}
made a systematic use of this method --- by exploiting the theory of stochastic
flows--- to study the solution of nonlinear first order partial differential
equations.
\\
We suppose first that $a^{0}=b^{0}=0$ and $f=g=0$ and denote by $\phi_{t}(x)$ the flow associated to the stochastic differential equation
\[ dx(t)= a(t,x(t))\circ dw(t)+b(x(t))dt. \]
Then if the initial condition is $C^{1}$, Kunita \cite{K1} showed that the  equation (${\cal E}_{S}$) has a unique global solution $u(t,x)$ in a strong sense, namely:
\[
u(t,x)=u_{0}(x)+\int_{0}^{1}a^{1}(\t,x)(\partial u/\partial x)\circ
dw(\t) +\int_{0}^{t}b^{1}(\t,x)(\t,x)(\partial u/\partial x)d\t,
\]
and furthermore, $u(t,x)=u_{0}(\phi^{-1}_{t}(x))$. Funaki showed that the last expression
gives a solution in a weak sense (for a similar equation that (${\cal E}_{S}$) with
boundary conditions).
\begin{prop}
Let $u_{0}\in H^{s}$. Then $u(t, \cdot)=u_{0}(\phi^{-1}_{t}(x))$ is the unique
solution to (${\cal E}$), in the case where $f=g=0$ and $a^{0}=b^{0}=0$.
\end{prop}
{\it Proof.} It suffices to approximate $u_{0}$ by a sequence $u^{n}_{0}\in
C^{1}(\R^{d})$ with $|u^{n}_{0}-u_{0}|_{s}\rightarrow 0$ as $n\rightarrow 
\infty$. Then the solution to the equation (${\cal E}_{S}$) with the initial data
$u^{n}_{0}$ is given by $u^{n}(t,x)=u_{0}^{n}(\phi^{-1}_{t}(x))$. Now using the
energy estimates of (\ref{energy1}), we get $E\sup_{t\in I}|u^{n}(t)-u(t)|_{s}^{2}
\rightarrow 0 $, and the proposition follows from the fact that
  $E\sup_{t\in I}|u^{n}(t)-u_{0}(\phi_{t}^{-1} (\cdot))
|_{s}^{2}\rightarrow 0 $ by Lebesgue's theorem. $\Box$
\\
\noindent
In the case where $a^{0}\neq 0$, $f\neq 0$ (and still $b^{0}=0,g=0$ for simplicity)
 the solution to (${\cal E}_{S}$) has the following expression:
\begin{align}
\label{kunita}
u(t,x)=& \{ u_{0}(\phi^{-1}_{t}(x))+ \\ \nonumber
  &+\int_{0}^{t}f(\t,\phi^{-1}_{\t,t} (\cdot))
    \exp(\int_{\t}^{t}a^{1}(r,\phi^{-1}_{r,t} (\cdot))\circ \hat{d}w(r))
     \circ \hat{d}w(\t)\}\\ \nonumber
     &\times\exp(\int_{0}^{t}a^{0}(\t,\phi^{-1}_{\t,t} (\cdot))\circ \hat{d}w(\t),
\end{align}
where $\int X(t)\circ \hat{d}w(t)$ denotes the backward Stratonovich integral,
taken in $H^{s-2}$. This expression follows from the same argument as in the 
above proof and the results of Kunita \cite{K1}.  
\\[12pt]
\noindent
{\bf (c) Remark on the use of the Fisk-Stratonovich differential.}
We want to show that the use of the Fisk-Stratonovich differential in Eq. ({\cal E})
is essential for obtaining the existence result. Let us consider a similar 
equation in which we use the It\^o differential: $ du(t)=a(x,D)u(t)dw(t)$.
Then we will have: 
\begin{align*}
 d\left< u(t),u(t) \PrR_{s} =& \left< u(t), (a(x,D)+a^{*}(x,D))u(t) \PrR_{s} dw(t)\\
   &+\frac{1}{2}\left< (a^{2}(x,D) +a^{* 2}(x,D))u(t),u(t) \PrR_{s} dt,
\end{align*}
and we can not obtain an energy estimate as in $\S$ 2.3 because the operator $a^{2}(x,D)+a^{*2}(x,D)$ is unbounded.

\noindent
In the case of scalar equations solved with the method of stochastic 
characteristics, the Fisk-Stratonovich notation is also essential for it 
allows the use of the same arguments as in the deterministic case.
\\[12pt]
\noindent
{\bf (d) Remark on the propagation speed.}
It is well known that the solution to deterministic symmetric systems of the 
form $ du(t)=a_{t}(x,D)u(t)dt,\; u(0)=u_{0}$ has a finite 
propagation speed, i.e. there is a constant $C>0$ such that if $u_{0}$
vanishes on $\{x: |x|>R\}$ then $u(t)$ will vanish on $\{x:|x|> R+Ct \}$. In 
the stochastic case, things are different: some coefficients are white
noises and thus 'unbounded'. Let us consider the simple
equation: 
\[ du(t)= \frac{\partial u}{\partial x}(x,t)\circ dw(t). \] 
Its solution is $u(t,x)=u_{0}(x+w(t))$. But since $\sup_{w}|w(t)|=
+\infty$ a.e. when $t>0$, we see that we can not have a finite propagation speed or a finite domain of dependence (\cite{Ta}) in
this case.

\subsection{Comparison with other existence results for SPDEs}

Hyperbolic SPDEs have been studied through several models. To cite only few examples, the stochastic wave equation in a space-time white noise setting is
one of the basic models, see Walsh \cite{Wa}; see also the interesting cases considered in Hajek \cite{Hajek} and \cite{Gaveau}. Hyperbolic equations or systems subject to additive noises are studied, e.g. in Chow \cite{Chow}, Dalang et al. \cite{DQS}, Lototsky-Rosovsky \cite{LR}. Kim \cite{Kim} considered a system of the form $\partial_{t}u+\sum_{i}A_{i}(t,x,u)\partial_{x_{i}}u=\sum_{i}f_{i}(u)dw_{i}$ where 
the $A_{i}$ are symmetric matrices and the $f_{i}$ are mappings that satisfy a Lipshitz condition on Sobolev spaces; this is a non linear model close to the one we study 
in this paper, but the last mentioned condition excludes the cases where the $f_{i}(u)$ are first order operators.

Ascanelli and S\"u\ss \, \cite{AS} studied a model of linear {\it scalar } hyperbolic SPDEs of the type: $L u(t, x) = \gamma(t,x) + \sigma(t,x)\dot{F}(t,x)$,
to which one has to give a sense in the framework of mild solution and stochastic integration with respect
to martingale measures; and in \cite{ACS}, Ascanelli, Coriasco and S\"u\ss \, considered a scalar stochastic hyperbolic equation of the type $L(t,x,\partial t,\partial x)u(t,x) = \gamma(t,x,u(t,x)) + \sigma(t,x,u(t, x)) \dot{\Xi}(t,x)$ with a space-time noise $\dot{\Xi}$. 

Finally, a class of hyperbolic-parabolic equations driven by standard Brownian motions received interest by many authors, e.g., in Lions, Perthame and Souganidis \cite{LPS1},  Bauzet, Vallet and Wittbold \cite{BVW} and Gess and Souganidis \cite{GS} where the deterministic entropy solution concepts are adapted to the stochastic case; in these cases, the equations are either scalar or the factor of the Brownian motion is a Lipshitz function of $u$.
\\
\noindent
These models can be included in one of the two main approaches to SPDEs: the first one considers SPDEs as stochastic evolution equations driven by Brownian motion in a Hilbert space \cite{DZ}, \cite{KR}, \cite{Pa}, \cite{Ro}. The second one considers partial differential equations perturbed by a space-time white noise (cf., e.g., Walsh \cite{Wa}). We refer to \cite{DQS} for an account and comparison of these approaches and to \cite{Kotelenez2}, \cite{LPS2}, \cite{PR}, \cite{GM} for further informations and references.
\noindent
As we are in the Hilbert space framework, we briefly explain why the standard existence results are not comparable to those of section 2.3. These results are essentially obtained by two methods: the variational method (see Pardoux \cite{Pa}, Krylov-Roszovskii \cite{KR},
Roszovskii \cite{Ro}) and the semi-group method (see Da Prato-Zabczyk \cite{DZ} and the references 
given there), and they are mainly concerned with parabolic type SPDEs.
\\[12pt]
{\bf (a) The variational method.} It is a generalization of the 
variational approach to PDEs; The framework is the following: let $V$ 
be a separable Banach space which is (continuously and
densely) imbedded in a Hilbert space $H$: $V\subset H\equiv H'\subset V'$ 
and we denote by $\|.\|$ and $|.|$ the norms in $V$ and $H$ respectively.
Now consider some operators $A\in {\cal L}(V,V'), B_{i}\in {\cal L}(V,H),
i=1,...,n$ and the equation:
\begin{equation}
\label{var1}
du(t)= A u(t)dt +B_{i}u(t) dw^{i}(t)
\end{equation}
with $u(0)=u_{0}\in H$ and $w^{i}(t), t\in I, i=1,...,n$ are standards 
independent Brownian motions defined on a probability space $(\Omega,{\cal F},
{\cal F}_{t}, P)$. The case of non linear equations is also treated by this approach. The main assumption is the 
following coercivity condition: there exist $\lambda, \gamma > 0$ such that
for all $v \in V$:
\[ -2\left<  Av, v \PrR_{V',V}+\lambda|v|^{2}\leq \gamma\|v\|^{2}+\sum_{i=1}^{n}
    |B_{i}v|^{2}. \]   
Under this condition, the equation (\ref{var1}) has a unique solution in
$M^{2}(I,V)$, the set  of adapted V-valued processes $u(t), t\in I$ with
$E\int_{I}\|u(t)\|^{2}< \infty$. 
In order to apply this result to our situation, we choose $V=H^{s}(\R^{d}),
H=H^{s-1}(\R^{d})$. For simplicity we choose $s=1$, so that $V=H^{1},
 H=L^{2}$. Now, let us consider the equation $du(t)= a(x,D)u(t)\circ dw(t)+
b(x,D)u(t)dt, u(0)=u_{0}\in H$ where $a, b\in {\rm OPS}^{1}$. This equation 
can be written in the It\^o form:
\begin{equation}
\label{eq252}
 du(t)=(\frac{1}{2}a(x,D)a(x,D)+b(x,D))u(t)dt+a(x,D)u(t)dw(t).
\end{equation}
The coercivity condition for this equation would be: there exist $\lambda, \gamma
>0$ such that:
\[  -\left< (a^{2}+2 b)v,v \PrR_{H^{-1}, H^{1}}+\lambda|v|_{L^{2}}^{2}\geq \gamma \|v\|_{H^{1}}^{2}+|av|_{L^{2}}^{2}. \]
\noindent
If we consider the simplest case where $a(x,D)u=\alpha \partial u/\partial x$,
$b(x,D)u=\beta \partial u/\partial x$, then, this condition implies that
$\lambda|u|_{L^{2}}^{2}\geq \gamma \|u\|_{H^{1}}^{2}$ for all $u\in H^{1}(\R^{d})$, which is not possible.
\\[12pt]
\noindent
{\bf (b) The semi-group method}. We consider again the equation
\begin{equation}
\label{semi1}
du(t)=A u(t)dt +B_{i}u(t)dw^{i}(t),
\end{equation}
Here $A$ is assumed to be the infinitesimal generator of a $C^{0}$-semigroup
 $S(t)$ in a Hilbert space $H$ and $u(0)=u_{0}\in H$ with 
$E|u_{0}|^{2}< \infty$. In this approach we look generally for a 
mild solution to Eq. (\ref{semi1}) i.e. $u(t)$ satisfies 
\[ u(t)= S(t)u_{0}+\int_{0}^{t} S(\t)B_{i}u(\t)dw^{i}(\t) .\]
Different assumptions are used; the first one is to
suppose that the operators $B_{i}$ are bounded, in which case (\ref{semi1})
has a unique mild solution. In the second one, the operators $B_{i}$ are
 allowed to be unbounded but the semigroup $S(t)$ is assumed to be analytic; a third one is 
to assume some Lipshitz conditions on $A, B_{i}$. These conditions are, however, not fulfilled in the case of $({\cal E})$ or
the example of equation (\ref{eq252}).
	
%%%%%%%%%%%%%%%%%%%%%%%%%%%%%%%%%%%%%%%%%%%%%%%%%%%%%%%%%%%%%%%%%%%%%%%%

\section{Small perturbations}
\subsection{Introduction and preliminaries}
This section is devoted to the study of the small random
perturbations of linear hyperbolic systems. More
precisely, let $u^{\epsilon} (\cdot)$ be the solution to 
\begin{equation}
\nonumber
({\cal S}_{\epsilon}):\left\{ \begin{array}{l}
  \displaystyle
 du^{\epsilon}(t) = \sqrt{\epsilon} a_{t}(x,D)u^{\epsilon}(t)\circ dw(t)+ b_{t}(x,D)u^{\epsilon}(t)dt,\\
\displaystyle
  u^{\epsilon}(0) = u_{0}\in H^{s},
\end{array}\right.
\end{equation}
where $a_{t},b_{t}$ are smooth families of (matrices of) pseudodifferential
operators which satisfy the conditions of Theorem 2.1. We denote by $u (\cdot)$
the solution of $({\cal S}_{0})$, the corresponding deterministic system. We 
are then interested in the limiting behavior of $u^{\epsilon} (\cdot)$ as $\epsilon\rightarrow 0$.
In the finite-dimensional case, problems of this type have been studied by many 
authors, see, e.g., Freidlin-Wentzell \cite{FW}, Deushel-Stroock \cite{DS} for references. In the infinite dimensional case, similar problems
have been addressed mainly for stochastic parabolic equations under various 
conditions, see Daprato-Zabczyk\cite{DZ} for references to earlier works on the subject.
In \cite{Ch}, Chow considered a small perturbation problem for 
the SPDE: 
\[du^{\epsilon}(t)=(Au^{\epsilon}(t)+F(u^{\epsilon}(t)))dt +\sqrt{\epsilon}
\Sigma(u^{\epsilon}(t))dw(t,
\]
where $A$ satisfies a coercivity condition (see $\S$ 2.4) and $\Sigma$ is
assumed to be Lipschitz in some sense. In \cite{Pe}, Peszat considered
the same problem in the semi-group framework with a set
of technical conditions which are not satisfied in our case. For other techniques that may be used in this context see also 
\cite{BDM}, \cite{RZ} and the references therein.

The method we use here is an adaptation of that of Priouret \cite{Priouret} in the finite dimensional 
case who follows an idea of Azencott \cite{Az}.
First, we state the following proposition which shows the convergence in probability of $u^{\epsilon}$ to $u$ with respect to the norm $\sup_{t}|v(t)|_{s-2}$:
\begin{prop}
For each $\delta >0$ we have 
\begin{equation}
\label{eq411}
\lim_{\epsilon\rightarrow 0}{\rm Pr}(\sup_{t\leq T}|u^{\epsilon}(t)-u(t)|_{s-2}
        >\delta       )=0 .
\end{equation}
\end{prop}
{\it Proof.} The proof is similar to the finite
dimensional case (see \cite{FW}); we shall give it in order to explain 
the loss of two derivatives in (\ref{eq411}). First, we recall that for 
$\epsilon>0$ bounded ($\leq 1$ say) we have from the previous section
\begin{equation}
\label{eq412}
E\sup_{t\leq T}|u^{\epsilon}(t)|^{4}_{s}< K<\infty ,
\end{equation}
for some $K>0$ (in particular $\sup_{t\leq T}|u(t)|^{4}_{s}\leq K$). Now by a 
 simple calculation we have
\begin{align*}
 |u^{\epsilon}(t)-u(t)|^{2}_{s-2}=& \sqrt{\epsilon}\int_{0}^{t}\sigma_{t}(
   u^{\epsilon}(\t))dw(\t)+\int_{0}^{t}\beta_{t}(u^{\epsilon})(\t)d\t\\
   &+\int_{0}^{t}\left< B_{\t}(x,D)(u^{\epsilon}(\t)-u(\t)),u^{\epsilon}(\t)-u(\t)
      \PrR_{s-2}d\t,
\end{align*}
with
\[ \sigma_{t}(v)= \left< A_{t}(x,D)v,v \PrR_{s-2}-2Re\left< v,a^{*}_{t}(x,D)u(t) \PrR_{s-2},\]
\[ \beta_{t}(v)= \frac{\epsilon}{2}[\left< L_{t}(x,D)v,v \PrR_{s-2}-
 2Re\left< v,a^{*2}_{t}(x,D)u(t) \PrR_{s-2}.\]
From the boundedness of the $B_{t}$ and the Gronwall lemma, it follows that
\begin{align*}
 \sup_{t\leq T}|u^{\epsilon}(t)-u(t)|^{2}_{s-2} \leq& C(T)[\frac{\epsilon}{2}
 \int_{0}^{T}(|u^{\epsilon}(t)|_{s-2}^{2}+|u^{\epsilon}(t)|_{s-2}|u(t)|_{s})dt
      \\ 
&+\sqrt{\epsilon}\sup_{t\leq T}|\int_{0}^{t}\sigma_{t}( u^{\epsilon}(\t))
dw(\t)|.
\end{align*}
From (\ref{eq412}) and the last inequality we get
\[ E\sup_{t\leq T}|u^{\epsilon}(t)-u(t)|^{2}_{s-2}\leq C'(T)\sqrt{\epsilon}\]
($C'(T)$ is another constant).
Now Proposition 3.1 follows from the last inequality. $\Box$

The objective is to give the exact rate of convergence in 
(\ref{eq411}). It turns out that this rate is exponential w.r.t $\epsilon$.
More precisely the family of the laws of $u^{\epsilon}$ satisfies a 
large deviation principle as in the finite dimensional case. 
\\
\noindent
Let $E$ be a topological space endowed with a $\sigma$-field ${\cal B}$. 
 We assume here that $E$ is Polish ${\cal B}$ is its Borel 
$\sigma$-field (although many results in large deviation theory hold in a 
more general setting). A function $I: E\longrightarrow [0, +\infty]$ is said
to be a rate function if it is lower semi-continuous. If in addition the
level sets $\{ x\in E: I(x)\leq L\}, L\geq 0$ are compact, then $I$ is 
said to be a good rate function.

\begin{defi}
A family $P^{\epsilon}, \epsilon >0$ of probability measures on 
$(E,{\cal B})$ satisfies a large deviation principle (LDP) with a rate 
function $I$ if
\[-\inf_{{\rm int} A} I(x)\leq \lim\inf_{\epsilon\rightarrow 0}
 \epsilon\log P^{\epsilon}(A) \leq \lim\sup_{\epsilon\rightarrow 0}
  \epsilon\log P^{\epsilon}(A)\leq -\inf_{{\rm cl} A} I(x), \]
for all $A\in {\cal B}$.
\end{defi}

We shall use the following standard result of large deviation theory (contraction principle):
\begin{prop}
Let $(E_{1},d_{1})$ and $(E_{2},d_{2})$ be two metric spaces and 
$X^{1}_{\epsilon}, X^{2}_{\epsilon}$ be two families of random variables
with values in  $E_{1}$ and $E_{2}$ respectively. Assume that the family
of laws ${\rm Pr}(X^{1}_{\epsilon}\in  \cdot)$ satisfies a large deviation principle 
with a good rate function $I$ and that there is a map 
 $\Phi: E_{1}\cap \{ I<+\infty\}\longrightarrow E_{2}$  such that:\\
(i) For all $L>0$, $\Phi_{|\{I\leq L\}}$ is continuous.\\
(ii) For each $h\in E_{1}$ with $I(h)<+\infty$ and $\eta>0$ we have
\[ \lim_{\delta \rightarrow 0}\limsup_{\epsilon\rightarrow 0} \epsilon
\log {\rm Pr} (d_{2}(X^{2}_{\epsilon},\Phi(h))>\eta, d_{1}(X^{1}_{\epsilon},h)< \delta ) =-\infty. \]
Then the family $Pr(X^{2}_{\epsilon}\in  \cdot)$ satisfies a large deviation
 principle with the good rate function
\[ I^{'}(y)=\inf_{x}\{I(x): \Phi(x)=y\}. \]
 \end{prop}

\subsection{A large deviation principle}
Let $u^{\epsilon}(t)$ be the solution to (${\cal S}_{\epsilon}$). For notational
simplicity we shall drop the index $t$ in $a_{t}, b_{t},$ etc. 
In this paragraph $u_{0}\in H^{s}$ is fixed and we denote by
${\cal H}^{s}$ the space $C_{u_{0}}(I, H^{s}), I=[0,T]$ of continuous paths
in $H^{s}$ starting from $u_{0}$. It will be equipped with the norm $\|v\|_{s,\infty} = \sup_{t\in I} |v(t)|_{s}$
and the corresponding Borel $\sigma$-field. Finally, let $P^{\epsilon}$ be the law of
$u^{\epsilon} (\cdot)$ which is defined on ${\cal H}^{s}$. It is also defined on 
all ${\cal H}^{s'}$ with  $s'<s$. We can now state the main result of this
section:
\begin{thm}
The family $(P^{\epsilon})$ satisfies a large deviation principle in 
${\cal H}^{s-2}$ with the following good rate function
\[ I_{u}(\phi)=\inf\{ \frac{1}{2}\int_{0}^{T}|\h(t)|^{2}dt : \Psi(h)=\phi \},
\]
where $\Psi : C_{0}([0,T], \R)\longrightarrow {\cal H}^{s}$ is given by
\begin{equation}
\label{psih}
 \Psi(h)(t)= u_{0}+\int_{0}^{t}a(x,D)\Psi(h)(\t)\h(\t)d\t+\int_{0}^{t}
b(x,D)\Psi(h)(\t)d\t.
\end{equation}
\end{thm}
First let us observe that the equation satisfied by $\Psi(h)(\cdot)$ is a deterministic hyperbolic system and has a unique 
solution in $H^{s}$. The rest of this paragraph is devoted to the proof of Theorem 3.4. 

The theorem will be proved by applying the contraction principle (Proposition 
3.3) with $\Phi=\Psi$, and $X^{1}_{\epsilon}=\sqrt{\epsilon}w (\cdot)$. 
From the Schilder theorem  we
know that $\mu^{\epsilon}$, the law of $X^{1}_{\epsilon}$, satisfies a LDP with 
the good rate function $I_{w}(h)= (1/2)\int_{0}^{T}|\h(t)|^{2}dt $.
Hence it suffices to verify the conditions (i) and (ii) of Proposition 3.3.

\begin{lm}
The map $\Psi:X:=(C_{0}([0,T], \R)\cap \{I(w)< \infty\}, |.|_{\infty})
\longrightarrow {\cal H}^{s-1}$ is continuous.
\end{lm}
{\it Proof.} Let $n\geq 0$ and $h\in H^{1}$. Define the polygonal 
approximation of $h$ by 
\[ h_{n}(t)= h([t]_{n})+(t-[t]_{n})\frac{h([t]_{n}+T/n)-h([t]_{n})}
{T/n},  \]
where we have used the following notation: if 
$t\in [iT/n, (i+1)T/n[$ then we set $[t]_{n}=iT/n$ 
(i.e. $[t_{n}]=[nt/T]T/n$). Now let $\Psi^{n}(h) (\cdot)$ be the
solution to the following equation:
\[ \Psi_{n}(h)(t)=u_{0}+\int_{0}^{t}( a(x,D)\Psi_{n}(h)(\t)\h_{n}(\t)
                          +b(x,D)\Psi_{n}(h)(\t))d\t \]
The map $\Psi_{n}: X\longrightarrow {\cal H}^{s-1}$ is continuous since 
$ \Psi_{n}(h)$ depends only on $h(iT/n), i=1,..., n$.

Next, we shall prove that for each $L>0$ the sequence $\Psi_{n}(h) (\cdot)$ 
converges uniformly (w.r.t. $h$) on the set $X_{L}:=
(C_{0}([0,T], \R)\cap \{I_{w}< L\} )$, namely:
\begin{equation}
\label{cvuniforme}
\lim_{n\rightarrow +\infty}\sup_{h\in X_{L}}\|\Psi_{n}(h) (\cdot)-\Psi(h) (\cdot)\|_{s-1}
=0 .
\end{equation}
Let $h\in X_{L}$. We have:
\begin{align*}
\Psi(h)(t)-\Psi_{n}(h)(t)=& \int_{0}^{t}(a(x,D)\h_{n}(\t)+b(x,D))
(\Psi(h)(\t)-\Psi_{n}(h)(\t))d\t\\
&+\int_{0}^{t}a(x,D)\Psi(h)(\t)(\h(\t)-\h_{n}(t))d\t
\end{align*}
and if we set $q_{n}(t)=\Psi(h)(t)-\Psi_{n}(h)(t)$, we get
\begin{align*}
d\left< q_{n}(t),q_{n}(t) \PrR_{s-1}=& (\left< (A(x,D)+B(x,D))q_{n}(t),q_{n}(t) \PrR_{s-1}
\h_{n}(t)\\
 &+2Re\left< q_{n}(t), a(x,D)q_{n}(t) \PrR_{s-1}(\h(t)-\h_{n}(t)))dt 
\end{align*}
Since $q_{n}(t)$ is uniformly bounded in $H^{s}$, it follows that
\begin{align*}
\phi_{n}(t)\leq & C\int_{0}^{t}\{ \phi_{n}(\t)|\h_{n}(\t)|d\t +
       \sqrt{\phi_{n}(t)}|\h(\t)-\h_{n}(\t)|\}d\t \\
     \leq& C\int_{0}^{t}( \phi_{n}(\t)|\h_{n}(\t)|+|\h(\t)-\h_{n}(\t)|)d\t
               +C\int_{0}^{T}|\h(\t)-\h_{n}(\t)|d\t
\end{align*}
But $\int_{0}^{T}|\h_{n}(t)|^{2}dt\leq \int_{0}^{T}|\h(t)|^{2}dt\leq L$ (by 
convexity), hence from the last inequality and the Gronwall lemma
it follows that:
\[ \phi_{n}(t)\leq C(\int_{0}^{T}|\h(\t)-\h_{n}(\t)|^{2}d\t)^{1/2} \]
On the other hand:
\begin{align*}
\int_{0}^{T}|\h(\t)-\h_{n}(\t)|^{2}d\t=& \sum_{i=1}^{n}\int_{iT/n}^{
       (i+1)T/n}(\h(\t)^{2}-\frac{(h((i+1)T/n)-
         h(iT/n))^{2}}{(T/n)^{2}})d\t\\
    =& \int_{0}^{T}\h(t)^{2}dt -\sum_{i=1}^{n}\frac{(h((i+1)T/n)-
         h(iT/n))^{2}}{T/n}.
\end{align*}
It is well known that the r.h.s of the above equality tends to $0$ as 
$n\rightarrow \infty$ for $h$ absolutely continuous and with derivative in
$L^{2}$. This implies that $\sup_{t}\phi_{n}(t)\rightarrow 0$ and 
the convergence is uniform on $\{ h: I_{w}(h)\leq L\}$. $\Box$

We shall prove that $\Psi$ satisfies the condition (ii) of the proposition.
First we consider the condition (ii) in the case $h=0$ and for a different map in the following lemma which is proved in the appendix:

\begin{lm}
Let $h\in X$ be fixed and consider $v^{\epsilon}(t)$ the solution to the 
equation
\[ dv^{\epsilon}(t)=\sqrt{\epsilon}a(x,D)v^{\epsilon}(t)\circ dw(t)+(b(x,D)+
a(x,D)\h(t))v^{\epsilon}(t)dt,\; v^{\epsilon}(0)=u_{0}.\] 
Then, we have for each $\eta >0$
\begin{equation}
\label{lm43}
 \limsup_{\delta\rightarrow 0}\epsilon\log {\rm Pr}(\sup_{t}
      |v^{\epsilon}(t)-\Psi(h)(t)|_{s-2}>\eta, |\sqrt{\epsilon}w|_{\infty}<
             \delta) =-\infty.
\end{equation}
($\Psi(h)$ is defined by (\ref{psih})).
\end{lm}
\noindent
Now, the passage from Lemma 3.2 to the condition (ii) of Proposition 3.3
can be done as in the finite-dimensional case by the following lemma; for completeness we give its
proof.
\begin{lm}
For each $h\in X$ and $\eta>0$, we have:
\begin{equation}
\label{lm44}
\lim_{\delta\rightarrow 0}\limsup_{\epsilon\rightarrow 0}\epsilon\log
  {\rm Pr}(\sup_{t}|u^{\epsilon}(t)-\Psi(h)(t)|_{s-2}\geq \eta,
     \sqrt{\epsilon}|w-h|_{\infty}< \delta )=-\infty,
\end{equation}
where $u^{\epsilon}$ is the solution to (${\cal E}_{\epsilon}$).
\end{lm}
{\it Proof.} For $\epsilon >0$ let $w^{\epsilon}(t)=w(t)-h(t)/\sqrt{\epsilon}, t\in [0,T]$. By the Girsanov theorem, $w^{\epsilon}$ is a Brownian motion
under the probability $Q^{\epsilon}$ with:
\[ \frac{d Q^{\epsilon}}{dP}=\exp(\int_{0}^{T}(\frac{\h(t)}{\sqrt{\epsilon}}
  dw(t)-\frac{1}{2}\int_{0}^{T}\frac{\h^{2}(t)}{\epsilon}dt ) .\]
Under the probability $Q^{\epsilon}$, $u^{\epsilon}$ satisfies 
\[ du^{\epsilon}(t)=\sqrt{\epsilon}a(x,D)u^{\epsilon}(t)\circ dw(t)+
     (\h(t)a(x,D)+b(x,D))u^{\epsilon}(t) dt\]
By Lemma 3.2 we have proved that (ii) holds for the above equation when $h=0$
i.e. 
\begin{equation}
\label{lm441}
\lim_{\delta\rightarrow 0}\limsup_{\epsilon\rightarrow 0}\epsilon\log
  Q^{\epsilon}(F(\epsilon,\eta,\delta))=-\infty,
\end{equation}
with $F(\epsilon,\eta,\delta))=\{\sup_{t}|u^{\epsilon}(t)-\Psi(h)(t)|_{s-2}
\geq \eta, \sqrt{\epsilon}|w^{\epsilon}|_{\infty}< \delta \}$. To prove the
lemma we have to show that $\lim_{\delta\rightarrow 0}\limsup_{\epsilon
\rightarrow 0}\epsilon\log P(F(\epsilon,\eta,\delta))=-\infty$. But
\[ P(F(\epsilon,\eta,\delta))=E^{Q^{\epsilon}}1_{F(\epsilon,\eta,\delta)}
  \frac{dP}{dQ^{\epsilon}}\leq
        (Q^{\epsilon}(F(\epsilon,\eta,\delta)))^{1/2}
   (E(\frac{dP}{dQ^{\epsilon}})^{2})^{1/2}.\]
Using the fact that $E^{Q^{\epsilon}}\exp(-\int_{0}^{T}2\h(t)/\sqrt{\epsilon}
  dw(t)-1/2\int_{0}^{T}4\h^{2}(t)/\epsilon dt)=1$, it follows that
\[ P(F(\epsilon,\eta,\delta))\leq Q^{\epsilon}(F(\epsilon,\eta,\delta))^{1/2}
   \exp\int_{0}^{T}\frac{\h^{2}(t)}{\epsilon}dt, \]
and
\[ \epsilon\log  P(F(\epsilon,\eta,\delta))\leq \frac{\epsilon}{2}\log 
  Q^{\epsilon}(F(\epsilon,\eta,\delta))+\frac{1}{2}\int_{0}^{T}\h^{2}(t)dt. \]
Now, (\ref{lm44}) follows from (\ref{lm441}) and the last inequality.
 $\Box$

%corrections 

%%%%%%%%%%%%%%%%%%%%%%%%%%%%%%%%%%%%%%%%%%%%%%%%%%%%%%%%%%%%%%%%%%%%%%%%

\section{Pathwise approximation and applications}
In this part, we consider the problem of pathwise approximation, also
called Wong-Zakai \cite{WZ} or Stroock-Varadhan \cite{SV} approximation in the case of SDEs. This aims at approximating
the stochastic solutions by the solutions of ordinary differential equations where the Brownian motion is regularized.
To be specific, let $(t_{i}^{n}),n\geq 1, 0\leq i\leq n$ be the subdivision of the interval
$[0,T]$ with $t^{n}_{i}=iT/n$ and consider the equation:
\[
({\cal E}_{n}):
\left\{ \begin{array}{c}
\displaystyle 
\frac{\partial u}{\partial t}=a_{t}(x,D)u(t) \w^{n}(t)+b_{t}( x,D)u +f(t,x)\w^{n}(t) +g(t,x)\\
u(0, \cdot)=u_{0} (\cdot),
\end{array}\right.
\]
where $w^{n}(t)$ is the polygonal approximation of the Brownian motion given
by: 
\[ w^{n}(t)=w([t]) +(t-[t])\frac{ \Delta w(t)}{\Delta t}. \]
We use here the following notation: if $t\in [t^{n}_{i}, t^{n}_{i+1}[$,
then we set 

\[[t]=t^{n}_{i}, \Delta w(t)= w(t_{i+1}^{n})-w(t_{i}^{n}) \]
and
\[\Delta t= t^{n}_{i+1}-t_{i}^{n}=T/n.\]

For each $w$, the (deterministic) equation (${\cal E}_{n}$) has a unique solution in $C^{0}(I, H^{s})$ which we denote
by $u^{n}$. Then we are interested in the convergence of $(u^{n})$ to
$u$. As we have mentioned in the introduction, this kind of approximation 
has been extensively studied in the case of SDEs, see, e.g., \cite{Doss}. In the case of SPDEs,
Gy\"ongy (see e.g., \cite{Gy}) studied this problem in an abstract 
variational framework which concerns parabolic SPDEs while Twardovska 
\cite{Tw} obtained other results by using the semi-group method. 
Br\'zezniak and Flandoli addressed this problem in the case of
{\it scalar} parabolic (possibly degenerate) SPDEs with an
application to scalar first order equation; in fact they use a representation
of the solutions to these equation via a Feynman-Kac type formula 
which reduces the problem to proving the approximation for a stochastic
flow of an associated SDE. In \cite{Roth}, Roth considered scalar stochastic hyperbolic equations for which 
a pathwise approximation is used
together with finite difference scheme in order to approximate the solutions.
More recently, Hairer and Pardoux \cite{HP} studied the case of non linear parabolic SPDEs driven by a space-time white noise 
by addressing the issue of the Stratonovich integration for space-time Brownian motion, 
and in \cite{Yastrzhembskiy}, Yastrzhembskiy considered the SPDE:  
$du(t, x) = [a^{ij} (t, x)D_{ij} u(t, x) +f(u, t, x)] dt +\sum_{k=1}^{m}g_{k}(u(t, x))dw^{k}(t)$, 
for which he proved a pathwise approximation result and a Stroock-Varadhan type support theorem in a suitable path space. 

These results do not seem to be applicable for the hyperbolic systems (${\cal E}$). Instead, we observe that
 the approximation is valid if the operator are bounded, like in the
finite-dimensional case. In the general case, we approximate the operators
by a family of bounded operators and we prove a uniform estimate (Lemma 4.1).
\noindent
In section 4.2 we prove a support type theorem for the SPDE (${\cal E}$) which extends the Stroock-Varadhan support theorem for SDEs (see \cite{SV},
\cite{IW}) to the infinite dimensional case of hyperbolic systems (${\cal E}$), and in section 4.3, we mention an application of the pathwise approximation to 
the random semigroup associated to Eq. (${\cal E}$). The results of this section, especially Proposition \ref{RandomOperator}, are used in \cite{AA} to extend the H\"ormander propagation of singularities theorem for the stochastic hyperbolic equations considered in this paper.

\subsection{Wong-Zakai type approximation}
 The purpose of this section is to prove the following
\begin{thm}
Under the assumptions of Theorem 2.1 we have:
\[ \lim_{n\rightarrow \infty} E\sup_{t\in I}|u^{n}(t)-u(t)|_{s-2}^{2}=0. \]
\end{thm}
For the sake of simplicity, the proof will be done in the following case:
$b=0, f=g=0$ and we shall assume that the operator $a_{t}(x,D)$ does
not depend on $t$. It will appear that the proof is valid for the setting of section 2. 

\noindent
Let $J_{\epsilon}, \epsilon \in ]0,1]$ be a Friedrichs mollifier and 
consider the solutions $u^{\epsilon}$ and $u^{\epsilon,n}$ to the equations
\[ ({\cal E}_{\epsilon}) : \; \; d u^{\epsilon}(t)= a(x,D) J_{\epsilon}u^{\epsilon}(t)\circ dw(t) \]
\[ ({\cal E}_{\epsilon,n}) :\; \; d u^{\epsilon,n}(t)= a(x,D)J_{\epsilon} u^{\epsilon,n}(t)
\w^{n}(t)dt \]
with the initial conditions $u^{\epsilon}(0)=u^{\epsilon,n}(0)=u_{0}$. Let
$\y= u^{\epsilon,n}-u^{\epsilon}$. The proof of Theorem 4.1 will be based on 
the following theorem and lemma.
\begin{thm}
For each fixed $\epsilon > 0$ we have
\[ \lim_{n\rightarrow \infty} E\sup_{t\in I}|u^{\epsilon,n}(t)-u^{\epsilon}(t
)|_{s}^{2} =0 .\]
\end{thm}
{\it Proof.} As the operators $a(x,D) J_{\epsilon}$ are bounded for each $\epsilon$ fixed, this theorem is proved in the same way
as the corresponding result of Nakao-Yamato \cite{NY} in the case of stochastic differential equations. 
\begin{lm}
There exist two functions $\alpha(\epsilon),\beta(n)$ with $\alpha(\epsilon)
\rightarrow 0$ as $\epsilon\rightarrow 0$ and $\beta(n)\rightarrow 0$ as
$n\rightarrow \infty$, such that
\[ E\sup_{t\in I}|u^{\epsilon,n}(t)-u^{n}(t)|^{2}_{s-2}\leq 
\alpha(\epsilon)(1+\beta(n)). \]
\end{lm}
The proof of this lemma is made in the Appendix.
\\
\noindent
{\it Proof of theorem 4.1.} Let $\delta >0$. By Lemma 4.1, there exists 
$\epsilon_{1}>0$  such that $E\sup_{t\in I}|u^{\epsilon_{1},n}(t)-u^{n}(t)|_{s-2}^{
2}< \delta $ for all $n\geq 1$. On the other hand, by Lemmas \ref{lm-Cauchy} and \ref{regul-limit} we can choose 
$\epsilon_{1}$ such that $E\sup_{t\in I}|u^{\epsilon_{1}}(t)-u(t)|_{s-2}^{2} < 
\delta$. Now Theorem 4.2 with $\epsilon=\epsilon_{1}$ implies that there
is $N\geq 1$ such that $E\sup_{t\in I}|u^{\epsilon_{1},n}(t)-u^{\epsilon_{1}}(t)|_{s}^{2}<\delta$ for all $n\geq N$. 
Summarizing,
\begin{eqnarray*}
E\sup_{t\in I}|u^{n}(t)-u(t)|_{s-2}^{2} &\leq& 3 (E\sup_{t\in I}|u^{n}(t)- 
                u^{\epsilon_{1},n}(t)|_{s-2}^{2} +E\sup_{t\in I}|u^{\epsilon_{1},n}
                  (t)-u^{\epsilon_{1}}(t)|_{s}^{2}\\
               &+& E\sup_{t\in I}|u^{\epsilon_{1}}(t)-u(t)|_{s-2}^{2})\\
               &\leq& 9\delta
\end{eqnarray*}
for $n\geq N$, which completes the proof of Theorem 4.1. $\Box$

%%%%%%%%%%%%%%%%%%%%%%
\subsection{Application to a support theorem}
In this section we apply the previous pathwise approximation to prove a support theorem
for the equation (${\cal E}$). Several extensions of this theorem to infinite dimensional settings have been
carried out. In \cite{Aida}, Aida proved a support theorem for diffusions in a Hilbert space where some of the complications 
of the infinite dimension are pointed out, and Nakayama \cite{Nakayama} proved a support theorem for the mild solution to equations 
of the type $dX(t) = AX(t)dt + b(X(t))dt + \sigma(X(t))dB(t)$, where $W$ is a cylindrical Brownian motion, $A$ is the infinitesimal 
generator of a ($C^{0}$)-semigroup $(S(t), t\geq 0)$ of bounded linear operators on $H$, and $b, \sigma$ are bounded and Lipshitz 
(which, as mentioned before, is different from the case we consider: for the system (${\cal E}$), the factor $\sigma(X(t))$ is a linear unbounded operator)  .

In the case of parabolic SPDEs similar results have been obtained
by Gy\"ongy \cite{Gy2}; the same problem has
been addressed for space-time white noise driven SPDEs of hyperbolic 
type by A. Millet and M. Sanz-Sol\'e \cite{MS1}.

As in Section 3, let ${\cal H}_{u_{0}}^{s}$ be the space of continuous
maps $I\longrightarrow (H^{s}(\R^{d}))^{d'}$ endowed with the norm
$\sup_{t\in I} |v|_{s}$ and the corresponding Borel $\sigma$-field. We 
know that the law $P_{u_{0}}$ of the solution to Equation (${\cal E}$) defines
a probability measure on $ {\cal H}_{u_{0}}^{s-2}$. In this section we 
describe the topological support of $P_{u_{0}}$, that is, the smallest closed 
subset $A$ of ${\cal H}_{u_{0}}^{s-2}$ such that $P_{u_{0}}(A)=1$.
\noindent
Let $H^{\infty}$ and $H^{\infty}_{p}$ be the sets of $\phi : I\longrightarrow \R$ which are respectively infinitely differentiable 
and piecewise infinitely differentiable, with $\phi(0)=0$. To each such function we associate the 
solution $(v(t,\phi), t\in I)$ to the following first order system
\begin{equation}
\label{support1}
d v(t)= a_{t}(x,D)v(t)\dot{\phi}(t)dt+b(x,D)v(t)dt,\; \; v(0)=u_{0} 
\end{equation}
As in the finite-dimensional case, we define
\[{\cal S}^{u_{0}}=\{v(t,\phi): \; \phi\in H^{\infty} \},\; \;
  {\cal S}^{u_{0}}_{p}=\{v(t,\phi): \; \phi\in H^{\infty}_{p}\},\]
In this section we shall assume that the family $a_{t}(x,D)$ satisfies the additional condition:
\\ \\
(iv')the operators $M(t):= L(t)a(t)+a^{*}(t)L(t)$ form a bounded family in 
${\rm OPS}^{0}$.
\\
Then we have the following result which extends the Stroock and Varadhan support theorem (\cite{SV}, \cite{IW}) to
the infinite dimensional case of Eq. (${\cal E}$):
\begin{thm}
The support ${\cal S}(P_{u_{0}})$ of $P_{u_{0}}$ is equal to $\overline{{\cal S}}^{u_{0}}= \overline{{\cal S}}^{u_{0}}_{p}$ where the closure is taken in 
${\cal H}^{s-2}_{u_{0}}$.
\end{thm}
{\it Proof.} First, it is easy to see that  $\overline{{\cal S}}^{u_{0}}= 
\overline{{\cal S}}^{u_{0}}_{p}$. Next, for $n\geq 1$, let $P^{n}_{u_{0}}$ be the law
of the solution $u^{n}$ to $({\cal E}_{n})$. Then $Q^{n}(\overline{{\cal S}}^{u_{0}}
)= Q^{n}(\overline{{\cal S}}^{u_{0}}_{p})=1$ for all $n$. By Theorem 4.1 it 
follows that $P^{n}_{u_{0}}\rightarrow P_{u_{0}}$ weakly, which implies
that $P_{u_{0}}(\overline{{\cal S}}^{u_{0}})\geq \limsup P^{n}_{u_{0}}(
\overline{{\cal S}}^{u_{0}})=1 $. Hence ${\cal S}(P_{u_{0}})\subset \overline{{\cal S}}^{
u_{0}}$.
\noindent
The inclusion $\overline{{\cal S}}^{u_{0}}\subset {\cal S}(P_{u_{0}})$ is a 
consequence of the following 
\begin{thm}
Let $\phi\in H^{\infty}$ and $\eta>0$. Then under the additional assumption (iv') we have
\begin{equation}
\label{support2}
\lim_{\delta\rightarrow 0 }{\rm Pr}(\sup_{t\in I}|u(t,w)-u(t,\phi)|_{s-2}>\eta
  ||w-\phi|_{\infty}<\delta)=0 .
\end{equation}
\end{thm}
This theorem is proved in the Appendix; the proof consists in reducing the theorem to the case of a stochastic evolution equation
with bounded operators which could be treated as in the finite dimensional case. 

%%%%%%%%%%%%%%%%%%%%%%
\subsection{Application to the random evolution operator}
One of the main applications of Wong-Zakai type approximations in the case
of a stochastic differential equation is the construction of its stochastic 
flow of diffeomorphisms. In this paragraph we give a similar application
which will be used later in the study of the singularities of Eq. (${\cal E}$). For notational simplicity, we 
still consider the case where $b=0,f=g=0$. Let $U_{f}(t',t)\phi$ denotes the 
solution to the forward equation:
\[ ({\cal E}_{F}):\; \; u(t)=\phi +\int_{t'}^{t}a_{\t}(x,D) u(\t)\circ dw(\t),\;
0\leq t'\leq t\leq T, \]
where $\phi\in H^{s}$ (we could also assume that $\phi$ is random and 
${\cal F}_{t'}$-measurable). Now, let us consider the backward equation
\[ ({\cal E}_{B}):\; \; u(t')=\phi -\int_{t'}^{t}a_{\t}(x,D) u(\t)\circ \hat{d}
w(\t),\; 0\leq t'\leq t\leq T, \]
where $\phi \in H^{s}$. For a fixed
$t\in ]0,T]$ we denote by ${\cal F}_{t',t}$ the $\sigma$-field 
$\sigma(w(\t)-w(\t'),t'\leq \t'\leq \t\leq t)$.

\begin{prop}\label{RandomOperator}
(i) The equation $({\cal E}_{B})$ has a unique solution $(u(t'))_{0\leq t'
\leq t}$ which is ${\cal F}_{t',t}$-adapted. We denote it by $U_{b}(t,t')\phi$.
\\
(ii) We have $U_{b}(t,t')U_{f}(t',t)=U_{b}(t,t')U_{f}(t',t)=Id$, a.e.
\end{prop}
{\it Proof.} (i) is proved exactly as in the case of forward equations. Also, the Wong-Zakai approximations holds for backward equations.\\
(ii)Let us denote by $U^{n}_{f}(t',t)\phi$ and $U^{n}_{b}(t',t)\phi$ the 
solutions to the following equations
\[ u(t)=\phi +\int_{t'}^{t}a_{\t}(x,D)u(\t)\w^{n}(\t)d\t, \]
\[ u(t')=\phi -\int_{t'}^{t}a_{\t}(x,D)u(\t)\w^{n}(\t)d\t. \]
Then we have for all $\phi\in H^{s}$:
\[ U^{n}_{b}(t,t')U^{n}_{f}(t',t)\phi=U^{n}_{f}(t',t)U^{n}_{b}(t,t')\phi=\phi.\]
Now the assertion (ii) follows from the approximation theorem of this section.
$\Box$

%%%%%%%%%%%%%%%%%%%%%%%%%%%%%%%%%%%%%%%%%%%%%%%%%%%%%%%%%%%%%%%%%%%%%%%%%%%%%%%%%%%%%%%%%%%%

%%%%%%%%%%%%%%%%%%%%%%%%%%%%%%%%%%%%%%%%%%%%%%%%%%%%%%%%%%%%%%%%%%%%%%%%%%%%%%%%%%%%%%%%%%%%%

\section{On the regularity of laws of the solutions}
In this section we consider the equation:
\begin{equation}
\label{reg1}
u(t)=u_{0}+\int_{0}^{t}a_{\t}(x,D)u(\t)\circ dw(\t)+\int_{0}^{t}b_{\t}(x,D)
  u(\t)d\t ,
\end{equation}
in the scalar case, where $u_{0}\in H^{s}(\R^{d})$ and $a_{t}(x,D), b_{t}(x,D)$
are smooth and bounded families of pseudodifferential operators in ${\rm OPS}^{1}$
such that their principal symbols are imaginary. By Theorem 2.1,
there is a unique solution to (\ref{reg1}) in $M^{2}(I,H^{s})$. We shall 
assume that $s > d/2$ so that $x\mapsto u(t,x)$ is continuous for each $t$. 
Then we are interested in the regularity of the law of the random variable
$u(t,x)$ for a given $(t,x)$. Similar problems of regularity of laws have
been studied for other classes of stochastic partial differential equations
of parabolic type in the case of one parameter driven white noise such as 
the Zakai equation of nonlinear filtering. On the other hand the same problems
have been addressed for parabolic and hyperbolic equations in the 
case of space-time noise (and one space dimension). See Nualart \cite{Nua} for
references.

We recall now some definitions and notations of the Malliavin calculus. Let
$X$ be a Hilbert space. We denote by $S(X)$ the set of ``simple'' $X$-valued
random variables $F$ of the form
\[ F(w)= f(w(t_{1}), \ldots, w(t_{n})),\; \; 0\geq t_{1}< \ldots < t_{n}\leq T,
\]
where $f: \R^{n}\longrightarrow X$ is a function which, together with all its partial derivatives, has a polynomial growth. We denote by $H$ the 
Cameron-Martin space i.e. $H:=\{ h\in H^{1}(\R): h(0)=0\}$. Given $h\in H$
and $F\in X$, the derivative of $F$ in the direction $h$ is defined by
\[ D_{h}F(w):=\frac{d}{d\epsilon}|_{\epsilon=0}F(w+\epsilon h)= 
 \sum_{i=1}^{n}\partial_{i} f(w) h(t_{i}).\]
The gradient of an $X$-valued random variables is the operator $D: S(X)
\longrightarrow L^{2}(\Omega\times [0,1]\times X)$ defined by
\[ D_{\theta}F:= \sum_{i=1}^{n}\partial_{i}f(w) 1_{\theta\leq t_{i}}, \]
so that 
\[ D_{h}F = \int_{0}^{1} D_{\theta}F.\h(\theta) d\theta, \; \; h\in H .\]
The operator $D$ is closable in $L^{p}(\Omega, X), p\geq 1 $ and 
 $\D^{1,2}$ will designate the domain of its closure in 
$L^{2}(\Omega\times X)$ (still denoted by $D$). We recall the following 
criterion of absolute continuity of laws (see \cite{Nua}):
\begin{prop}
Let $F$ be a real valued random variable in $\D^{1,2}(\R)$. Assume that
\begin{equation}
\label{reg2}
  \|DF\|_{L^{2}([0,1])} > 0 \; \; a.e. 
\end{equation}
then the law of $F$ is absolutely continuous with respect to the Lebesgue
measure.
\end{prop}
Now we return to Equation (\ref{reg1}). We first state the following

\begin{prop}
Let $u$ be the solution to Eq. (\ref{reg1}). Then for each $t\in I$ we have
$u(t)\in \D^{1,2}(H^{s-2})$ with $D_{\theta} u(t)=0 $ if $\theta > t$ and
\begin{equation}
\label{reg3}
D_{\theta} u(t)= a_{\theta}(x,D)u(\theta)+\int_{\theta}^{t}a_{\theta}(x,D)D_{
\theta} u(\t)\circ dw(\t)+\int_{\theta}^{t}b_{\theta}(x,D)D_{\theta} u(\t)\circ d\t
\end{equation}
In other words,
\begin{equation}
\label{reg4}
D_{\theta}u(t)= U(\theta, t)a_{\theta}(x,D) u(\theta),
\end{equation}
where $U(\theta, t)$ is the stochastic evolution semigroup associated
te Eq. (\ref{reg1}).
\end{prop}
{\it Proof.} For the sake of simplification, we suppose that $b\equiv 0$. Let 
$\epsilon>0$ and $u^{\epsilon}$ be the solution to
\[ u^{\epsilon}(t)=u_{0}+\int_{0}^{t}a_{\t}(x,D)J_{\epsilon}u^{\epsilon}(\t)\circ dw(\t).\]
Since the operators $a_{t}(x,D)J_{\epsilon}$ are bounded in $H^{s}$, it can be shown 
that $u^{\epsilon}(t)\in \D^{1,2}(H^{s})$ for $t\in I$ and
\[ D_{\theta} u^{\epsilon}(t)= a_{\theta}(x,D)J_{\epsilon}u^{\epsilon}(
\theta)+\int_{\theta}^{t}a_{\theta}(x,D)D_{\theta}J_{\epsilon} u(\t)
\circ dw(\t). \]
The proof of these is exactly the same as for finite dimensional SDEs, see 
e.g. \cite{Nua}.
Next we shall show that $D u^{\epsilon}(t)$ is a Cauchy family in $L^{2}(\Omega
\times I\times H^{s-2})$. Before doing this, we remark that
\begin{equation}
\label{reg5}
\sup_{ 0<\epsilon< 1}\sup_{\theta\in I} E |D_{\theta}
u^{\epsilon}(t)|_{s-1}^{4}<+\infty,
\end{equation}
Indeed, for $\theta$ fixed we have
\[ d |D_{\theta} u^{\epsilon}(t)|_{s-1}^{2}=\left< A_{\epsilon}(t)D_{\theta}  u^{\epsilon}(t),  D_{\theta} u^{\epsilon}(t) \PrR_{s-1}dw(t)
+\frac{1}{2}\left<    L_{\epsilon}(t) D_{\theta} u^{\epsilon}(t),D_{\theta} u^{\epsilon}(t) \PrR_{s-1}, \]
where we have used the notation of section 2. Hence
\[ E |D_{\theta} u^{\epsilon}(t)|_{s-1}^{4}\leq C(E |u^{\epsilon}(\theta)
    |_{s}^{4} +\int_{0}^{\t} E|D_{\theta}u^{\epsilon}(\t)|_{s-1}^{4}d\t \]
which implies (\ref{reg5}). Now to show that $D u^{\epsilon}(t)$ is a Cauchy family, we have to
estimate $\v^{\theta}(t):= D_{\theta} u^{\epsilon}(t)-D_{\theta} u^{\epsilon'}
(t) $ which satisfies
\[ d\v^{\theta}(t)=a_{\t}(x,D) J_{\epsilon}\v(t)\circ dw(t)+a_{\t}(x,D)( 
  J_{\epsilon}- J_{\epsilon'})D_{\theta}u^{\epsilon'}(t)\circ dw(t). \]
To this end, we use (\ref{reg5}) and the same calculations as in 2.3.2 (b). 
We omit the details. Now we have: $u(t)\in \D^{1,2}(H^{s-2})$ and 
$D u(t)$ is the limit of $D u^{\epsilon}(t)$ in $L^{2}(\Omega\times I\times 
H^{s-2})$. The fact that $D_{\theta} u(t)$ satisfies (\ref{reg3}) for 
$0\leq \theta\leq t$ can be proved easily by showing that 
$E \sup_{t\in I}|D_{\theta} u^{\epsilon}(t)-v_{\theta}(t)|_{s-2}^{2}
\rightarrow 0$ as $\epsilon\rightarrow 0$; here $v_{\theta}(t)$ is the 
solution to $dv_{\theta}(t)=a_{t}(x,D)v_{\theta}(t)\circ dw(t)$ for $t\geq 
\theta$ and $v_{\theta}(\theta)=a_{\theta}(x,D) u(\theta)$. $\Box$

\begin{Rq}{\rm 
We also see that for $h\in H$, $D_{h} u(t)$ satisfies
\[ D_{h}u(t)=\int_{0}^{t}a_{\t}(x,D)D_{h} u(\t)\circ dw(\t)+
        \int_{0}^{t}a_{\t}(x,D)u(\t)\h(\t)d\t 
\]
From this equation, we deduce an equivalent form of
(\ref{reg4}) (``Duhamel's principle'')
\[ D_{h}u(t)=\int_{0}^{t}U(\t,t)a_{\t}(x,D)u(\t)\h(\t)d\t.\]}
\end{Rq}
We return now to the absolute continuity of the law of $u(t,x)$ for a
given $(t,x)$. We assume from now on that $s-1>d/2$. This implies that
$u(t,x)\in \D^{1,2}(\R)$. To see this, we note that if $s>d/2$ then
$u^{\epsilon}(t,x)\in D^{1,2}(\R)$ and (\ref{reg3}) can be written 
`pointwise'. This is done as in the finite dimensional case; indeed
it suffices to estimate the pointwise norms by Sobolev norms. Now 
assuming that $s-1> d/2$, we get from  (\ref{reg5}) that
\begin{equation}
\label{reg6}
\sup_{\theta}\sup_{\epsilon} E |D_{\theta} u^{\epsilon}(t,x)|^{4}
   \leq C\sup_{\theta}\sup_{\epsilon} E |D_{\theta} u^{\epsilon}(t)|^{4}_{s}
  < \infty.
\end{equation}
Since $u^{\epsilon}(t,x)\rightarrow u(t,x)$ in $L^{2}(\Omega)$, we deduce
from (\ref{reg6}) that $u(t)\in \D^{1,2}(\R)$, by using e.g. Lemma 1.2.3 in
Nualart \cite{Nua}. Furthermore, from (\ref{reg3}) we deduce that
\[ D_{\theta} u(t,x)= (U(\theta,t)a_{\theta}(x,D)u(\theta))(x),\; \;
 0\leq \theta\leq t . \]
Using Proposition 5.1 we see that a sufficient condition for the law of
$u(t,x)$ to be absolutely continuous w.r.t. the Lebesgue measure is 
\begin{equation}
\label{reg7}
\int_{0}^{t} |(U(\theta,t)a _{\theta}(x,D)u(\theta))(x)|^{2}d\theta >0, \; a.e.
\end{equation}
Since we know that $U(\theta,t): H^{s}\longrightarrow H^{s}$ is continuous for 
almost all $w$ we deduce that $\theta \mapsto U(\theta,t)a _{\theta}(x,D)
 u(\theta)$ is continuous (w.r.t. the norm $|.|_{H^{s-1}}$) and
then $\theta\mapsto (U(\theta,t)a _{\theta}(x,D)u(\theta))(x)$ is 
continuous a.s. Hence a sufficient condition to have (\ref{reg7}) is
\[ |U(t,t)a_{t}(x,D)u(t)(x)|:=|a_{t}(x,D)u(t)(x)|>0 \; a.s.\]
or
\[ |U(0,t)a_{0}(x,D)u_{0}(x)|>0 \; a.s. \]
To go further, let us see the particular case of differential operators:
$a_{t}(x,D):= a^{i}(t,x)\partial/\partial x^{i}$, $b(x,D):=b^{i}(t,x)
\partial/\partial x^{i} $ for which $U(0,t) v(x)$ is given by Eq. 
(\ref{kunita}). This implies that a sufficient condition for (\ref{reg7}) to hold is that
\[ a_{0}(x,D) u_{0}(x) \neq 0 \; \; \mbox{for all} \; x.\]
In fact, to deduce this last condition we have used the strict positivity of
the semi-group $U(0,t)$ in the case of differential operator i.e.
 $U(0,t)\phi (x)>0$ for all $x$ whenever $\phi$ is continuous and $\phi(x)>0$
for all $x$. Now the remaining question is whether the semigroup $U(0,t)$ is
``strictly positive'' when $a(x,D)$ is a pseudodifferential operator.

%%%%%%%%%%%%%%%%%%%%%%%%%%%%%%%%%%%%%%%%%%%%%%%%%%%%%%%%%%%%%%%%%%%%%%%%

\section{ Appendix: Proofs of some technical lemmas}
In order to simplify the proofs and notations, we suppose that the operators $a_{t}(x,D)$ and $b_{t}(x,D)$ do not depend
on the time variable $t$, and they will be denoted by $a, b$; this will be therefore the case
for the operators $A_{t}, B_{t}$ and $L_{t}$ which will be denoted by $A, B$ and $L$.

\subsection{Proofs of the lemmas related to small perturbations}
{\it Proof of Lemma 3.2.}: 
The proof of this lemma will use the following:
\begin{lm}\label{LemmaA1}
Let $Z_{i}(t), i=1,2,3, t\in [0,T]$ be three adapted processes such that
\begin{equation}
\nonumber
\left\{ \begin{array}{l}
  \displaystyle
dZ_{1}(t) = \epsilon Z_{2}(t)\circ dw(t)+Z_{3}(t)dt, \\
\displaystyle
 |Z_{2}(t)| \leq A_{c} Z_{1}(t) \; {\rm with} \; A_{c} > 0,  \\   
\displaystyle
\int_{0}^{T} [\frac{|Z_{3}(t)|}{Z_{1}(\t)} + \frac{|Z_{2}(t)|}{Z_{1}(t)^{2}} + \frac{|\left< Z_{2},w\PrR_{t}|}{Z_{1}(t)} ]dt  \leq K \; a.e.
\end{array}\right.
\end{equation}
Then, assuming $Z_{1}(0)=1$, we have
\[ {\rm Pr}(\sup_{t\in [0,T]}|Z_{1}(t)| \geq M)\leq \exp[-\frac{(\log M - K)^{2}}{2T\epsilon^{2}A_{c}^{2}}].\]
\end{lm}
{\it Proof.} Let $\gamma >0$. By the stochastic calculus rules and It\^o formula we have:
\begin{align*}
d \log(Z_{1}(t)+\gamma)=& \frac{1}{Z_{1}(t)+\gamma} \circ dZ_{1}(t)\\
                     =&  \epsilon \frac{Z_{2}(t)}{Z_{1}(t)+\gamma}  \circ dw(t)+\frac{Z_{3}(t)}{Z_{1}(t)+\gamma}dt \\
                     =& \epsilon \frac{Z_{2}(t)}{Z_{1}(t)+\gamma}  dw(t) -\epsilon^{2} \frac{Z_{2}^{2}(t)}{2(Z_{1}(t)+\gamma)^{2}}dt 			
							+ \epsilon\frac{\left< Z_{2},w\PrR_{t}}{2(Z_{1}(t)+\gamma)}+\frac{Z_{3}(t)}{Z_{1}(t)+\gamma}dt
\end{align*}
Then, using the assumptions of the lemma we get for $\epsilon \leq 1$:
\begin{equation}
\label{EqA1}
| \log(Z_{1}(t)+\gamma)|\leq |\log(1+\gamma)|+K+\sup_{t\leq T}\epsilon |\int_{0}^{t}\frac{Z_{2}(\t)}{ Z_{1}(\t)+\gamma}dw(\t) |.
\end{equation}
Now the lemma follows from the exponential inequality for martingales; we recall a particular case of this inequality that will be also used in the proofs of some lemmas below: if $M_{t}$ is a martingale such that $\left<M\right>_{t} \leq c t, \; \forall t \in [0, T]$ for some constant $c$, then for $a\geq 0$: 
\begin{equation}
\label{ExpIneqMartingale}
{\rm Pr} (\sup_{t\leq T} M_{t} \geq a t) \leq e^{-a^{2}T/2c}
\end{equation} 
See, e.g. Revuz-Yor \cite{RevuzYor} (Exercise 3.16, p. 145). In our case we take:
\[ M_{t} = \epsilon |\int_{0}^{t}\frac{Z_{2}(\t)}{ Z_{1}(\t)+\gamma}dw(\t)| \]
and we have:
\begin{align*} 
\left<M\right>_{t} =& \epsilon^{2} \int_{0}^{t}\frac{Z_{2}^{2}(\t)}{ (Z_{1}(\t)+\gamma)^{2}}d\t \\
                   \leq & \epsilon^{2}  A_{c}^{2} t
\end{align*}
By letting $\gamma \longrightarrow 0$ in (\ref{EqA1}) we get:
\begin{align*} 
{\rm Pr}(\sup_{t\in [0,T]}|Z_{1}(t)| \geq M)  \leq & \; {\rm Pr}(\sup_{t\in [0,T]} M_{t} \geq \frac{(\log M - K)}{T} T) \\
                                    \leq &  \exp[-\frac{(\log M - K)^{2}}{T^{2}}\times\frac{T}{2\epsilon^{2}  A_{c}^{2}}]
																		=\exp[-\frac{(\log M - K)^{2}}{2T\epsilon^{2}A_{c}^{2}}]							
\end{align*}
where we have applied (\ref{ExpIneqMartingale}) with $c=\epsilon^{2}  A_{c}^{2}$ and $a=(\log M - K)/T$. $\Box$
\\
\noindent
Now we turn to the proof of Lemma 3.2. Let us denote $q_{\epsilon}(t)=v^{\epsilon}(t)-
\Psi(h)(t)$ and:
\[ F(\epsilon,\eta,\delta)=\{\sup_{t}
      |v^{\epsilon}(t)-\Psi(h)(t)|_{s-2}>\eta, |\sqrt{\epsilon}w|_{\infty}<
             \delta) \}.
\]
First, by a standard localization argument, $v^{\epsilon}(t)$ may be assumed bounded in $H^{s-2}$. 
Indeed, let $\t^{\epsilon}$ be the stopping time defined by:
\[ \t^{\epsilon}=\inf\{ t: |v^{\epsilon}(t)-\Psi(h)(t)|_{s-2}\geq \eta\}
     \wedge T\]
then we have:
\[ \{ \sup_{t\leq T}|v^{\epsilon}(t)-\Psi(h)(t)|_{s-2}\geq \eta\}=\{
              \sup_{t\leq \t^{\epsilon}}|v^{\epsilon}(t)-\Psi(h)(t)|_{s-2}
\geq \eta\}.\]
But if $t\leq \t^{\epsilon}$ then $|v^{\epsilon}(t)|_{s-2}\leq \sup_{t\leq T}|
    \Psi(h)(t)|_{s-2}+\eta=: M$.
Hence we have: 
\begin{equation}
\label{lm431}
 F(\epsilon,\eta,\delta)=F\cap \{t\leq \t^{\epsilon}\}\subset F\cap 
      \{\sup_{t\leq T}|v^{\epsilon}(t)|_{s}\leq M \}. 
\end{equation}
Next, we have:
\begin{align*}
v^{\epsilon}(t)-\Psi(h)(t)=& \int_{0}^{t}\sqrt{\epsilon}a(x,D)u^{\epsilon}
          (\t)\circ dw(\t) \\
  &+\int_{0}^{t}(a(x,D)\h(\t)+b(x,D))(u^{\epsilon}(\t)-
              \Psi(h)(\t))\h(\t)d\t,
\end{align*}
which yields:
\begin{align*}
\left< q_{\epsilon}(t),q_{\epsilon}(t) \PrR_{s-2}=& \sqrt{\epsilon}\int_{0}^{t}
       \sigma(v^{\epsilon}(\t))dw(\t)+\int_{0}^{t}\beta(v^{\epsilon}(\t))d\t\\
    &+\int_{0}^{t}\left< q_{\epsilon}(\t), (\h(\t)A(x,D)+B(x,D))q_{\epsilon}(\t))  \PrR_{s-2} d\t,
\end{align*}
where $q_{\epsilon}(t)=v^{\epsilon}(t)-\Psi(h)(t)$ and
\begin{align*}
\sigma(v_{\epsilon}(t)) =& \left< v^{\epsilon}(t),A(x,D)v^{\epsilon}(t) \PrR_{s-2}
        -2Re\left< v^{\epsilon}(t),a^{*}(x,D)\Psi(h)(t) \PrR_{s-2},\\
\beta(v_{\epsilon}(t)) =& \left< L(x,D)v^{\epsilon}(\t), v^{\epsilon}(\t) \PrR_{s-2}
      +\left< v^{\epsilon}(\t), a^{*2}(x,D)\Psi(h)(\t) \PrR_{s-2}
\end{align*}
Using the boundedness of $A(x,D),B(x,D)$ and $\int_{0}^{T}\h^{2}(t)dt$ it 
follows by the Gronwall lemma that:
\begin{align*}
 q_{\epsilon}(t) \leq & C\sup_{\theta\in [0,T]}|\sqrt{\epsilon}\int_{0}^{
\theta}\sigma(v_{\epsilon}(\t))dw(\t)| \\
&+\epsilon\int_{0}^{T}|\beta(v_{\epsilon}(\t))|d\t.
\end{align*} 
Then, according to (\ref{lm431}), we have $F(\epsilon,\eta,\delta)=F_{1}
\cup F_{2}$ with 
\[
 F_{1}=\{\sup_{\theta\leq T}|\sqrt{\epsilon}|\int_{0}^{\theta\wedge 
  \t^{\epsilon}}\sigma(v_{\epsilon}(\t))dw(\t)| \geq \eta^{2}/C ,
\sqrt{\epsilon}|w|_{\infty}<\delta\}, 
\]
\[
F_{2}=\{\epsilon\int_{0}^{T}|\beta(v_{\epsilon}(\t))|d\t \geq \eta^{2}/C,
 \sup_{t\leq T}|v^{\epsilon}(t)|_{s}\leq M, \sqrt{\epsilon}|w|_{\infty}<\delta
 \}. 
\]
By the boundedness of $L(x,D)$ and $\sup_{t}|\Psi(h)(t)|_{s}$, we have
$F_{2}(\epsilon,\eta,\delta)=\emptyset$ for $\epsilon\leq \epsilon_{0}$
with $\delta_{0}$ sufficiently small.
Now, for $n\geq 1$ we set:
$v^{n,\epsilon}(t)=v^{\epsilon}([t]_{n})$ and we have for $\gamma>0$:
\[ F_{1}(\epsilon,\delta,\eta)\subset A(\epsilon,\gamma,n)\cup
    B(\epsilon,\eta, \gamma,n)\cup C(\epsilon,\eta, \delta,n) \]
with:
\begin{align*}
 A(\epsilon,\gamma,n)=& \{\sup_{t\leq \t^{\epsilon}}|v^{\epsilon}(t)-
     v^{\epsilon,n}(t)|_{s-1} >\gamma \},\\
B(\epsilon,\eta, \gamma,n)=& \{\sup_{t\leq \t^{\epsilon}}|v^{\epsilon}(t)-
     v^{\epsilon,n}(t)|_{s-2} \leq\gamma,\\ 
      & \sup_{t\leq T}|\sqrt{\epsilon}\int_{0}^{t}
 (\sigma(v^{\epsilon}(\t))-\sigma(v^{\epsilon,n}(\t)))dw(\t)|\geq \eta^{2}/2C
      \}, \\
C(\epsilon,\eta, \delta,n)=& \{\sup_{t\leq T}|\sqrt{\epsilon}\int_{0}^{t}
           \sigma(v^{\epsilon,n}(\t)))dw(\t)|\geq \eta^{2}/2C, \sqrt{\epsilon}
                  \sup_{t}|w(t)|\leq \delta \}.
\end{align*}
First, observe that $\sigma$ is uniformly Lipshitz (constant $k$) on 
$\{ y\in H^{s-2}: |y|_{s-2}\leq M  \}$; hence if 
$|v(t)-v^{\epsilon,n}(t)|_{s}\leq \gamma$
and $|v(t)|_{s}\leq M$ then $\sqrt{\epsilon}|\sigma(v^{\epsilon}(t))-
\sigma(v^{\epsilon,n}(t))|\leq k\gamma\sqrt\epsilon$. Then, by the exponential
inequality of martingales (\ref{ExpIneqMartingale}), we have: 
\begin{equation}
\label{lm432}
 {\rm Pr}(B(\epsilon,\eta,\gamma,n))\leq 2\exp(
     -\frac{\eta^{4}}{8C^{2}k^{2}\gamma^{2}\epsilon}).
\end{equation}
Next, we turn to estimate ${\rm Pr}(A(\epsilon,\gamma,n))$. We have
\[ v^{\epsilon}(t)-v^{\epsilon,n}(t)=\sqrt{\epsilon}\int_{[t]_{n}}^{t}
     a(x,D) v^{\epsilon}(\t)\circ dw(\t)+\int_{[t]_{n}}^{t}(a(x,D)
 v^{\epsilon}(\t)\h(t)+b(x,D)v^{\epsilon}(\t))d\t. \]
hence:
\[
\left< v^{\epsilon}(t)-v^{\epsilon,n}(t),v^{\epsilon}(t)-v^{\epsilon,n}(t) \PrR_{s-2}=
\int_{[t]_{n}}^{t}Y(v^{\epsilon})
     (\t)dw(\t) +\int_{[t]_{n}}^{t}Z(v^{\epsilon}) (\t)d\t, \]
with:
\begin{align*}
Y(v^{\epsilon})(\t)=& \sqrt{\epsilon}
   \left< A(x,D)v^{\epsilon}(\t),v^{\epsilon}(\t) \PrR_{s-2}-2Re\left< a(x,D)v^{\epsilon}(\t)
   \PrR_{s-2},\\
Z(v^{\epsilon})(\t)  =& \frac{1}{2}\left< L(x,D)v^{\epsilon}(\t),v^{\epsilon}(\t) \PrR_{s-2}
       -2Re \left<  a^{2}(x,D)v^{\epsilon}(\t),v^{\epsilon}(\t) \PrR_{s-2}\\
    +&2Re\left<  v^{\epsilon}(\t)-v^{\epsilon,n}(\t), a(x,D)v^{\epsilon}(\t)\h(\t)+
         b(x,D) v^{\epsilon}(\t) \PrR_{s-2}. 
\end{align*}
On the other hand we have:
\[ {\rm Pr}(A(\epsilon,\gamma,n)\leq {\rm Pr}(A(\epsilon,\gamma,n), \sup_{t\leq T}|v^{\epsilon}(t)|_{s}\leq M)+ {\rm Pr}(\sup_{t\leq T}
|v^{\epsilon}(t)|_{s}\geq M) .\]
Now let $Z_{1}(t)=|v^{\epsilon}(t)|_{s}^{2}$. By the It\^o formula it follows
 that: 
\begin{equation} 
\label{EqZ1}
dZ_{1}(t)=\sqrt{\epsilon} Z_{2}(t)\circ dw(t)+Z_{3}(t)dt
\end{equation}
where $Z_{2}, Z_{3}$ are given by:
\begin{align*}
Z_{2}(t) & = \left< A(x,D)v^{\epsilon}(t),v^{\epsilon}(t) \PrR_{s} \\
 Z_{3}(t) &=\frac{1}{2}\left< L(x,D)v^{\epsilon}(\t),v^{\epsilon}(\t) \PrR_{s} + \left< (B(x,D)+A(x,D)\h(t)) v^{\epsilon}(\t),v^{\epsilon}(\t) \PrR_{s}.
\end{align*}
\noindent
By the boundedness of the operators $ A, B, L$, we see that  $Z_{1}, Z_{2}$ and $Z_{3}$ satisfy the assumptions of Lemma \ref{LemmaA1} (which we use with $\sqrt{\epsilon}$ in the Eq. (\ref{EqZ1}) instead of the corresponding $\epsilon$ of the above-mentioned lemma; we also use the fact that $\int_{0}^{T} |\h(t)| dt \leq \int_{0}^{T}(1+ \h(t)^{2})dt < \infty $). This implies that there is a constant K such that:
\begin{equation}
\label{veps}
 {\rm Pr}(\sup_{t\leq T}|v^{\epsilon}(t)|_{s}\geq M) \leq
    \exp(-\frac{(\log M-K)^{2}}{4T A_{c}\epsilon}).
\end{equation}
Where $A_{c}$ is a bound of the operator $A=a+a^{*}$. Next we have
\[ Pr(A(\epsilon,\gamma,n),\sup_{t\leq T}|v^{\epsilon}(t)|_{s}\leq M ) \leq S_{1}+S_{2}, \]
where:
\begin{eqnarray*}
S_{1} &=& \sum_{i=1}^{n-1}Pr(\sqrt{\epsilon}\sup_{t\in
   [iT/n,(i+1)T/n]}|\int_{i T/n}^{t\wedge \t^{\epsilon}}Y(v^{\epsilon}(\t))
       dw(\t) |> \gamma^{2}/4, \sup_{t\leq T}|v^{\epsilon}(t)|_{s}\leq M),\\
S_{2} &= &\sum_{i=1}^{n-1}Pr(\sup_{t\in
   [iT/n,(i+1)T/n]}|\int_{i T/n}^{t\wedge \t^{\epsilon}}|Z(v^{\epsilon})(\t)|
    d\t>\gamma^{2}/4, \sup_{t\leq T}|v^{\epsilon}(t)|_{s}\leq M).
\end{eqnarray*}
Noting that for $t\leq \t^{\epsilon}$ we have that $|Z((v^{\epsilon})(t))|
 \leq CM(1+|\h(t)|)$ for some constant $C$, we get:
\[ S_{2}\leq n{\rm Pr}(\int_{iT/n}^{(i+1)T/n}CM(1+|h(t)|)dt>\gamma^{2}/4)=0,\]
for $n\geq n_{0}$ sufficiently large. Also, concerning $S_{1}$, we have 
$|Y(v^{\epsilon}(t))|\leq C M$ and by the exponential inequality (\ref{ExpIneqMartingale}), we get
\[ S_{1}\leq \sum_{i=1}^{n-1} 2\exp(-\frac{\gamma^{4}}{32\epsilon C^{2}M^{2}
     n^{-1}}),
 \]
hence, for $n\geq n_{0}$ we have:
\begin{equation}
\label{lm433}
 {\rm Pr}(A(\epsilon,\gamma,n)\leq 2n\exp(-\frac{n\gamma^{4}}{32\epsilon 
   C^{2}M^{2}}\leq\exp(- \frac{n\gamma^{4}}{64\epsilon C^{2}M^{2}}),
\end{equation}
provided that $n_{0}$ is sufficiently large. 

Finally, as regards $C(\epsilon,\eta,\delta,n)$, noting that on $\{\sqrt{\epsilon}|w|_{\infty} \leq \delta \}$:
\[ |\sqrt{\epsilon}\int_{0}^{t}\sigma(v^{\epsilon,n}(\t)))dw(\t)|\leq
  \sqrt{\epsilon} |\sum_{i=1}^{n}\sigma(v(i\Delta_{n}))(w((i+1)\Delta_{n})-
            w(i\Delta_{n})|\leq 2C' Mn\delta,\]
we get: $C=\emptyset$ if $\delta < \eta^{2}C'/MCn$.
To summarize, let $R>0$. By (\ref{lm432}) there is a real $\gamma_{0}>0$ such that:
for all $n\geq 1$ we have
\[ {\rm Pr}(B(\epsilon,\eta,\gamma_{0},n))\leq \exp(\frac{-R}{\epsilon}).\]
By (\ref{veps}) there exists $M>0$ such that 
\[ {\rm Pr}( \sup_{t\leq T}|v^{\epsilon}(t)|^{2}_{s}\geq M)\leq \exp(-\frac{R}{
\epsilon}). \]
By (\ref{lm433}), ($\gamma_{0}, M$ being fixed) there exists $n_{1}\geq n_{0}$ 
such that
\[ {\rm Pr}(A(\epsilon,\gamma,n_{1},\sup_{t\leq T}|v^{\epsilon}(t)|^{2}_{s}
\leq M )\leq \exp(-\frac{R}{\epsilon}), \]
and if we choose $\delta\leq \delta_{0}:= \eta^{2}C'/CMn_{1}$ (so that
${\rm Pr}(C(\epsilon,\eta,\delta,n_{1}))=0$ we get
\[ {\rm Pr} (F(\epsilon,\delta,\eta))\leq 2\exp(\frac{-R}{\epsilon}),\]
i.e $\epsilon\log {\rm Pr} F(\epsilon,\delta,\eta)\leq -R+2\epsilon$, which
completes the proof of Lemma 3.2. $\Box$

\subsection{Proofs of the lemmas related to pathwise approximation and support theorem}
\subsubsection{Proof of Lemma 4.1}
First we state the following lemma which will be used in the proof.
\begin{lm}\label{LemmaA2}
Let $u^{n}(t), u^{\epsilon,n}(t)$ be the solutions to ${\cal E}_{n}, {\cal E}_{\epsilon,n}$ respectively, 
with the same initial value $u_{0}\in H^{s}$. Then
\[ E\sup_{t\in I} |u^{n}(t)|_{s}^{8}+ E\sup_{t\in I} |u^{\epsilon,n}(t)|_{s}^{8}\leq C, \]
where $C$ is a constant which depends only on $E|u_{0}|_{s}^{8}$ (and not
on $n, \epsilon$).
\end{lm}
{\it Proof.} We have:
\begin{align*}
\left< u^{n}(t), u^{n}(t) \PrR_{s}=& |u_{0}|_{s}^{2}+\int_{0}^{t}\left< A u^{n}(\t), u^{n}(\t) \PrR\w^{n}(\t)d\t\\
 =& |u_{0}|_{s}^{2}+\int_{0}^{t}\left< A u^{n}([\t]), u^{n}([\t]) \PrR\w^{n}([\t])d\t\\
  &+\int_{0}^{t}(\left< A u^{n}(\t), u^{n}([\t]) \PrR-\left< A u^{n}([\t]), u^{n}([\t]) \PrR)
\w^{n}([\t])d\t\\
=& |u_{0}|_{s}^{2}+\int_{0}^{t}\left< A u^{n}([\t]), u^{n}([\t]) \PrR dw(\t)\\
&+\int_{0}^{t}\left<  L u^{n}(c_{\t}), u^{n}(c_{\t}) \PrR(\w^{n}(\t))^{2}(\t-[\t])d\t,
 \end{align*}
where $c_{\t}\in ][\t],\t[$. Let $p^{n}(t)=\sup_{\t\leq t}|u^{n}(\t)|_{s}^{8}
$, then using the boundedness of $A,L$ and the martingale inequality, we
get:
\begin{align*}
E\sup_{t\leq t}|u^{n}(\t)|_{s}^{8}\leq &C\int_{0}^{t} 
E\sup_{\theta\leq t}|u^{n}(\theta)|_{s}^{8} d\t \\
+&CE(\int_{0}^{t}|u_{n}(c_{\t})|^{2}(\Delta w(\t))^{2}\frac{(\t-[\t])}{(
\Delta \t)^{2}}d\t)^{4}.
\end{align*}
The use of Schwarz's inequality for the last term does not permit to
conclude (via the Gronwall lemma). However we have:
\[ \mbox{Claim:}\;\; 
  E(\int_{0}^{t}|u_{n}(c_{\t})|^{2}(\Delta w(\t))^{2}\frac{(\t-[\t])}{(
\Delta \t)^{2}}d\t)^{4}\leq C \int_{0}^{t}\sup_{\theta\leq t}|u^{n}(\theta)|_{s}^{8} \frac{(\t-[\t])^{4}}{(\Delta \t)^{4}} d\t. \]
Using this and denoting $\psi^{n}(t)=E\sup_{t\leq t}|u^{n}(\t)|_{s}^{8}$ we
get
\[ \psi^{n}(t)\leq C\int_{0}^{t}\psi^{n}(\t)(1+\frac{(\t-[\t])^{4}}{(\Delta \t)^{4}})d\t. \]
Since: 
\[\int_{0}^{T}\frac{(\t-[\t])^{4}}{(\Delta \t)^{4}}d\t= \frac{T}{\Delta \t}
 \sum_{i=1}^{n}\int_{\t_{i}}^{\t_{i+1}}\frac{(\t-[\t])^{4}}{(\Delta \t)^{4}} d\t= \frac{T}{5},\]
 we get by the Gronwall lemma: 
\[ E\sup_{t\leq T}|u^{n}(t)|_{s}^{8}=\psi^{n}(T)\leq E|u_{0}|_{s}^{2}(
1+ e^{T^{2}/5}(1+T/5)).\]
{\it Proof of Claim.} We write for $\theta\in ][\t], [\t]^{+}[$:
\[\PrL u^{n}(\theta),u^{n}(\theta) \PrR =\PrL u^{n}([\theta]),u^{n}([\theta]) \PrR +
\int_{[\t]}^{\theta}\PrL Au^{n}(\lambda),u^{n}(\lambda) \PrR 
\frac{\Delta w(\t)}{\Delta \t} d\t .\]
and then:
\[   |u^{n}(\theta)|_{s}^{2}\leq |u^{n}([\t])|_{s}^{2} +C\int_{[\t]}^{\theta}
 |u^{n}(\lambda)|_{s}^{2}\frac{|\Delta w|}{\Delta \t}d\lambda. \]
By the Gronwall lemma this yields:
\[ |u^{n} (\theta)|_{s}^{2}\leq u^{n}([\t])|_{s}^{2}e^{\frac{A |\Delta w|}{\Delta \t} (\theta -[\t])}.\]
Now, using this inequality and the fact that the increment $\Delta w([\t])$ is
independent of the ${\cal F}_{[\t]}$ (and then on $u([\t])$) we get:
\begin{align*}
E(\int_{0}^{t} |u_{n}(c_{\t})|^{2}(\Delta w(\t))^{2}\frac{(\t-[\t])}{(
\Delta \t)^{2}}d\t)^{4} \leq & T^{3}\int_{0}^{t} d\t [E|u^{n}([\t])|^{8} \times \\
                     & \times E(e^{\frac{4 A |\Delta w|}{\Delta \t} (\t -[\t])}(\Delta w(\t))^{8} 
                                  \frac{(\t-[\t])^{4}}{(\Delta \t)^{8}}].
\end{align*}
But:
\begin{align*}
E e^{\frac{4 A |\Delta w|}{\Delta \t} (\t -[\t])}(\Delta w(\t))^{8} =&
\int \frac{x^{8}}{\sqrt{2\pi \Delta \t}}e^{\frac{4 A |x|}{
\Delta \t} (\t -[\t])}  e^{-\frac{x^{2}}{2\Delta \t}}dx\\
=& (\Delta \t)^{4}\int e^{\frac{4 A |x|}{\sqrt{\Delta \t}}
(\t-[\t])} e^{-x^{2}}{2} dx \\
\leq& (\Delta \t)^{4}\int e^{4 A\sqrt{T}|x| -
\frac{x^{2}}{2}} dx.
\end{align*}
Consequently,
\[ E(\int_{0}^{t}|u_{n}(c_{\t})|^{2}(\Delta w(\t))^{2}\frac{(\t-[\t])}{(
\Delta \t)^{2}}d\t)^{4}\leq C\int_{0}^{t}E\sup_{\theta\leq t}|u^{n}(\theta)|_{s}^{8} \frac{(\t-[\t])^{4}}{(\Delta \t)^{4}} d\t. \]
This proves the claim and the uniform boundedness of 
$E\sup_{t\in I}|u^{n}(t)|_{s}^{8}$. The proof is similar in the
 case of $u^{\epsilon,n}$. $\Box$

\noindent
{\it Proof of Lemma 4.1 .} We have, with $\y= u^{\epsilon,n}-u^{\epsilon}$: 
\[ d\y(t) = a \j \y(t) \w(t)dt +a(\j u^{n}-u^{n}) \w(t)dt \]
\begin{align*}
d\PrL\y(t),\y(t) \PrR =&  \PrL A_{\epsilon}\y(t),\y(t) \PrR \w(t)dt\\
 &+2\PrL\y(t),a (\j u^{n}(t)-u^{n}(t)) \PrR\w(t)dt,
\end{align*}
with $A_{\epsilon}=a\j+\j a^{*}$. We rewrite the above equation as
\begin{align*}
\PrL\y(t),\y(t) \PrR=& \int_{0}^{t}(\PrL  A_{\epsilon}\y([\t]),\y([\t]) \PrR\\
&+2\PrL \y([\t]), a (\j u^{n}([\t])-u^{n}([\t])) \PrR)\w(\t)d\t\\
 &+\int_{0}^{t}(\PrL  A_{\epsilon}\y(\t),\y(\t) \PrR
-\PrL  A_{\epsilon}\y([\t]),\y([\t]) \PrR)d\t\\
&+\int_{0}^{t}2(\PrL \y(\t), a (\j u^{n}(t)-u^{n}(\t)) \PrR\\
&-\PrL \y([\t]), a (\j u^{n}(t)-u^{n}([\t])) \PrR)\w(\t)d\t,
\end{align*}
Noting that for an adapted process $G$ we have $\int_{0}^{t}G([\t])\w(\t)d\t
=\int_{0}^{t}G([\t])dw(\t)$, we get
\begin{align*}
\PrL \y(t),\y(t) \PrR=& \int_{0}^{t}(\PrL  A_{\epsilon}\y([\t]),\y([\t]) \PrR\\
&+2\PrL \y([\t]), a (\j u^{n}([\t])-u^{n}([\t])) \PrR))dw(\t)\\
&+\int_{0}^{t}\{ \PrL L_{\epsilon} \y(c_{\t}),\y(c_{\t}) \PrR\\
&+\PrL \y(c_{\t}),K_{\epsilon}(\j u^{n}(c_{\t})-u^{n}(c_{\t})) \PrR \}(\w(\t))^{2}
(\t-[\t]) d\t \\
&+\int_{0}^{t}\{\PrL \y(c_{\t}'), a(\j a^{n}(c_{\t}')-a u^{n}(c_{\t}')) \PrR\\
&+\PrL \y(c_{\t}'), \j a^{*} a (\j u^{n}(c_{\t}')-u^{n}(c_{\t}')) \PrR\}(\w(\t))^{2}
(\t-[\t])d\t \\
&+\int_{0}^{t} \PrL  a (\j u^{n}(c_{\t}')-u^{n}(c_{\t}')), a(\j a^{n}(c_{\t}')
-a u^{n}(c_{\t}')) \PrR\\
&\times (\w(\t))^{2}(\t -[\t])d\t,
\end{align*}
where we have set $L_{\epsilon}=\j a^{*}A_{\epsilon} +A_{\epsilon}a\j,
K_{\epsilon}= A^{*}_{\epsilon}+A_{\epsilon}a$ and $c_{\t}, c_{\t}'\in 
][\t],\t[$. Let:

\[ \z(t)=\sup_{\t\leq t}|\y(\t)|_{s}^{4}. \]
\noindent
Then from the boundedness of the last two operators, martingale and Schwarz inequalities it
 follows that:
\begin{align*}
E \z(t)\leq &C E \int_{0}^{t}( \z(\t)+ \sqrt{\z(\t)}\psi^{1}_{\epsilon,n}(\t))
d\t \\
 &+CE \int_{0}^{t} (\z(\t)+\sqrt{\z(t)} \psi^{2}_{\epsilon,n}(\t)) (\w(\t))^{4}(\t-[\t])^{2}d\t \\
 &+E\int_{0}^{t}\psi^{3}_{\epsilon,n}(\t) (\w(\t))^{4}(\t-[\t])^{2}d\t\\
:=& I_{1}+I_{2}+I_{3},
\end{align*}
where 
\begin{align*}
\psi^{1}_{\epsilon,n}(t)=& \sup_{\t\leq t}|\j u^{n}(\t)-u^{n}(\t)|_{s+1},\\
\psi^{2}_{\epsilon,n}(t)=& \sup_{\t\leq t}| K_{\epsilon}
(\j u^{n}(t)-u^{n}(t))|_{s}^{2}+|a(\j a^{n}(t)-a u^{n}(t))|_{s}^{2}\\
&+| \j a^{*} a (\j u^{n}(t)-u^{n}(t))|_{s}^{2},\\
\psi^{3}_{\epsilon,n}(t)=&
  |a (\j u^{n}(t)-u^{n}(t))|^{2}_{s}|a(\j a u^{n}(t)-a u^{n}(t))|_{s}^{2}.
\end{align*}
Now, using Lemma 2.2 we get:
\[ E\sup_{\t\leq T} (\psi^{1}_{\epsilon,n}(t))^{2}\leq E\sup_{\t}|\j u^{n}(\t)-u^{n}(\t)|_{s+1}^{4}\leq \alpha_{1}(\epsilon) 
E\sup_{\t}|u^{n}(\t)|_{s+1}^{4},
\]
and similarly 
\[  E\sup_{\t\leq T} (\psi^{2}_{\epsilon,n}(t))^{2}\leq \alpha_{2}(\epsilon)
 E\sup_{\t}|u^{n}(\t)|_{s+2}^{8} , \]
\[ E\sup_{\t\leq T} (\psi^{3}_{\epsilon,n}(t))^{2}\leq \alpha_{3}(\epsilon)
 E\sup_{\t}|u^{n}(\t)|_{s+2}^{8} , \]
with $\alpha_{i}(\epsilon)\rightarrow 0$ as $\epsilon \rightarrow 0$, $i=1,2,3
$. Hence, denoting 
\[ \phi_{\epsilon,n}(t)= E\z(t), \]
it follows that
\[ I_{1}\leq  C\int_{0}^{t}(\phi_{\epsilon,n}(\t)+\sqrt{\phi_{\epsilon,n}(\t)}\alpha_{1}(\epsilon)) d\t. \]
For the term $I_{2}$, the direct use of the Schwarz inequality does not lead
to the good estimate; instead we proceed as follows:
\begin{align*}
I_{2} \leq & 16 \int_{0}^{t}(\phi_{\epsilon,n}(\t)+\sqrt{\phi_{\epsilon,n}(\t)} (E(\psi^{2}_{\epsilon,n
})^{2})^{1/2}) \frac{(\t-[\t])^{2}}{(\Delta \t)^{4}} d\t\\
 &+\int_{0}^{t}\{ [ (E\z(\t)^{2})^{1/2} +
  (E(\z(\t)))^{1/2} (E(\psi_{\epsilon,n}^{2}(
\t))^{2})^{1/2}] \\
&\times (E (\Delta w(\t))^{8} 1_{\Delta |w(\t)|\geq 2 \Delta 
\t})^{ 1/2} \frac{(\t-[\t])^{2}}{(\Delta \t)^{4}} \} d\t.
\end{align*}
Here we have simply used Schwarz's inequality and 
$E X (\Delta w(\t))^{4}\leq 16 (\Delta \t)^{2} E(X)
+ E(X 1_{(\Delta w(\t))^{2}\geq 4 \Delta \t})$ for a random variable $X$. We recall that $\w(t)=(\Delta w(t))/\Delta t$ with $\Delta t= T/n$. 
On the other hand we have:
\begin{align*}
E (\Delta w(t))^{8} 1_{(\Delta w)^{2}\geq 4 \Delta t}\leq & 
2\int_{2\sqrt{\Delta t}}^{+\infty} \frac{8}{\sqrt{2\pi \Delta t}}e^{-x^{2}/
2(\Delta t)} \\
 \leq & C e^{-1/\Delta t}.
\end{align*}
Now, using Lemma \ref{LemmaA2} which implies that $E(\z(T))^{2}, E\z(T)$ are uniformly bounded (w.r.t. $n$ and $\epsilon$), we get:
\begin{align*}
I_{2} \leq & C\int_{0}^{t}(\phi_{\epsilon,n}(\t)+\sqrt{\phi_{\epsilon,n}(\t)}\alpha_{2}(\epsilon)
   \frac{(\t-[\t])^{2}}{(\Delta \t)^{4}} d\t\\
 &+C\int_{0}^{t} \frac{-1/\Delta \t}{(\Delta \t)^{2}}\int_{0}^{t}\alpha_{2}
(\epsilon) (\t-[\t])^{2}d\t\\
\leq &C\int_{0}^{t}(\phi_{\epsilon,n}(\t)+\sqrt{\phi(\t)_{\epsilon,n}}\alpha_{2}(\epsilon)
   \frac{(\t-[\t])^{2}}{(\Delta \t)^{4}} d\t +\beta(\Delta \t),
\end{align*}
with $\beta(\Delta \t)\rightarrow 0$ as $\Delta t\rightarrow 0$. Finally, for
the term $I_{3}$ we have:
\begin{align*}
I_{3} \leq &\int_{0}^{t} (E(\psi^{3}_{\epsilon,n})^{2})^{1/2} (E(\Delta w)^{8})^{1/2} \frac{(\t-[\t])^{2}}{(\Delta \t)^{4}} \\
  \leq& C \alpha_{3}(\epsilon)\int_{0}^{t}\frac{(\t-[\t])^{2}}{(\Delta \t)^{2}} \leq C \alpha_{3}(\epsilon).
\end{align*}
Now, since we know that $\phi_{\epsilon,n}(t)$ is bounded, we can estimate the terms 
$\sqrt{\phi_{\epsilon,n}(t)}$ by a constant and we get an estimate of the form:
\[\phi_{\epsilon,n}(t)\leq C \int_{0}^{t}\phi_{\epsilon,n}(\t)(1+\frac{(\t-[\t])^{2}}{(\Delta \t)^{2}})
 d\t +\alpha(\epsilon)(1+\beta(\Delta \t)), \]
and the use of the Gronwall lemma completes the proof of Lemma 4.1. $\Box$

%%%%%%%%%  %%%%%%%%%%

\subsubsection{Proof of Theorem 4.4.}
Let $J_{n}$ be the Friedrichs mollifier $J_{1/n}$ and let us denote by $u(t,\phi)$ the solution to Eq. (\ref{support1}) and $u_{n}$ the 
solution to: 
\begin{equation}
\label{th44reduction}
du_{n}=a_{t}(x,D)J_{n}u_{n}\circ dw(t)+b_{t}(x,D)J_{n}u_{n}(t)dt,\; \; u_{n}(0)
=u_{0}
\end{equation}
To prove Theorem 4.4, we shall reduce its assertion to proving the same limit (\ref{support2}) for the solution $u_{n}$ to Eq. (\ref{th44reduction}), for which 
the operators $a_{t}(x,D)J_{n}$ and $b_{t}(x,D)J_{n}$ are bounded and this can be done as in finite dimension; this is the purpose of the lemma \ref{LemmaA5} below. 
The reduction to this case is done via the following two lemmas:

\begin{lm}\label{LemmaA3}
For $A>0$ sufficiently large we have:
\[ P_{1}(\delta, n):={\rm Pr}(\sup_{n}\sup_{t\in I}|u_{n}(t)|_{s}> A||w|_{\infty}
<\delta)\leq c \exp(-c'(\log A)^{2}) ,\]
( $c$ is independent of $\delta\in ]0,1]$). Also
\[  P_{1}'(\delta):={\rm Pr}(\sup_{t\in I}|u(t)|_{s}> A||w|_{\infty}
<\delta)\leq c \exp(-c'(\log A)^{2}) .\]
\end{lm}
\begin{lm}\label{LemmaA4}
There exists $N\geq 1$ such that: 
\[ \limsup_{\delta\rightarrow 0}{\Pr}(\sup_{t\in I}|u(t,w)-u_{N}(t,w)|_{s-2}>
    \eta, \sup_{n}\sup_{t}|u_{n}(t)|_{s}\leq A, \sup_{t\in I}|u(t)|_{s}\leq A
    | |w|_{\infty}<\delta) =0. \]
\end{lm}
\begin{lm}\label{LemmaA5}
Let $\eta >0$. Then for each $N>0$ fixed, we have:
\[ \lim_{\delta \rightarrow 0}{\rm Pr}(\sup_{t\in I}|u_{N}(t, w)-u(t,\phi)|_{s-2}>\eta,
  \sup_{n}\sup_{t}|u_{n}(t)|_{s}\leq A, \sup_{t\in I}|u(t)|_{s}\leq A
    | |w|_{\infty}<\delta)=0. \]
\end{lm}
For the sake of simplicity we will prove these lemmas in the case $b\equiv 0$.

\noindent
{\it Proof of Lemma \ref{LemmaA3}.}
Let $\gamma>0, n\geq 1$ and $Z_{n}(t):=|u_{n}(t)|_{s}^{2}$. Then by the It\^o
formula we have:
\begin{align*}
\log(Z_{n}(t)+\gamma)=& \log(|u_{0}|_{s}^{2}+\gamma)+\int_{0}^{t}\frac{\PrL A_{n}u_{n}(t), u_{n}(\t) \PrR _{s}}{Z_{n}(\t)+\gamma}\circ dw(\t) \\
   =& \log(|u_{0}|_{s}^{2}+\gamma)+ \frac{\PrL A_{n}u_{n}(\t), u_{n}(\t) \PrR_{s}w(t)}{Z_{n}(\t)+\gamma} \\
   &-\int_{0}^{t}\frac{\PrL L_{n}u_{n}(\t),u_{n}(\t) \PrR_{s}}{(Z_{n}(\t)+\gamma)}w(\t)\circ dw(\t)\\
    &+\int_{0}^{t}\frac{ \PrL A_{n}u_{n}(\t), u_{n}(\t) \PrR_{s}^{2}}{ (Z_{n}(\t)+\gamma)^{2}} w(\t)\circ dw(\t) \\
   =& \log(|u_{0}|_{s}^{2}+\gamma) + \frac{\PrL A_{n}u_{n}(\t), u_{n}(\t) \PrR_{s}w(\t)}{Z_{n}(\t)+\gamma} \\
   &-\int_{0}^{t}\frac{\PrL L_{n}u_{n}(\t),u_{n}(\t) \PrR_{s}}{(Z_{n}(\t)+\gamma)}w(\t) dw(\t) \\
    &+\int_{0}^{t}\frac{ \PrL A_{n}u_{n}(\t), u_{n}(\t) \PrR_{s}^{2}}{(Z_{n}(\t)+\gamma)^{2}} w(\t) dw(\t)\\
      &+\frac{1}{2} \int_{0}^{t}\{\frac{\PrL A_{n}u_{n}(\t), u_{n}(\t) \PrR_{s}^{2}}{(Z_{n}(\t)+\gamma)^{2}}
			             -\frac{\PrL L_{n}u_{n}(\t),u_{n}(\t) \PrR_{s}}{(Z_{n}(\t)+\gamma)} \\
      &+w(\t)[\frac{1}{(Z_{n}(t)+\gamma)^{2}}(\PrL L_{n}u_{n}(\t),u_{n}(\t) \PrR_{s}\PrL A_{n}(\t)u_{n}(\t),u_{n}(\t) \PrR_{s} \\
			&- \PrL M_{n}u_{n}(\t),u_{n}(\t) \PrR_{s}(Z_{n}(\t)+\gamma))\\
     &+\frac{2}{(Z_{n}(\t)+\gamma)^{2}}\PrL A_{n}u_{n}(\t), u_{n}(\t) \PrR_{s} \PrL L_{n}u_{n}(\t), u_{n}(\t) \PrR_{s} \\
          & - \frac{2}{(Z_{n}(\t)+\gamma)^{3}}\PrL A_{n}u_{n}(\t),u_{n}(\t) \PrR_{s}^{3} )]\} d\t\\
    :=& \log(|u_{0}|_{s}^{2}+\gamma) + k_{1,n}(t)w(t)+\int_{0}^{t}k_{2,n}(\t)w(\t)dw(\t) \\
		&+\int_{0}^{t}(k_{3,n} (\t)w(\t)+k_{4,n})d\t,
\end{align*}
where we have used the following notation: 
\begin{align*}
A_{n}&= a J_{n}+J_{n}a^{*},\\
   L_{n}&= A_{n}a J_{n}+J_{n}a^{*}A_{n}, \\ 
M_{n} &= L_{n} a J_{n}+J_{n}a^{*}L_{n}
\end{align*}
   
All these operators form a bounded family in ${\rm OPS}^{0}$ under the assumptions (iii)--(v), (iv'). Hence,
there is a constant $K$ such that $k_{i,n}(t)\leq K$ a.s. for all $n$ and $i=1,
...,4.$ Therefore, on the set $\{ |w|_{\infty}\leq \delta\}$, there 
is a constant $M>0$ such that:
\[ \log(Z_{n}(t))\leq M+\sup_{t\in I}|\int_{0}^{t}k_{2,n}(\t)w(\t)dw(\t)
  |  \]
Then for $A$ sufficiently large we get:
\[ {\rm Pr}(\sup_{n}\sup_{t\in I}\log(|u_{n}(t)|_{s}^{2})\geq (\log A),
     |w|_{\infty}<\delta )\leq c_{1}\exp -\frac{c_{2}(\log(A))^{2}}
    {\delta^{2}} \]
Now, we recall the following estimate of ${\rm Pr}(|w|_{\infty}<\delta)$ (see Ikeda-Watanabe \cite{IW}, lemma 8.1, p.519): there exist constants $c_{3}, c_{4}$ such that:
\begin{equation}
\label{iw1}
 {\rm Pr}(|w|_{\infty}<\delta)\sim c_{3}\exp-\frac{c_{4}}{\delta^{2}}. 
\end{equation}
It follows that for $A$ sufficiently large we have for some $c,c'>0$
\[ {\rm Pr}(\sup_{n}\sup_{t\in I}|u_{n}(t)|_{s}> A||w|_{\infty}
<\delta)\leq c \exp(-c'(\log A)^{2}) .\]
The second estimate of the lemma is proved in the same way as the first one. $\Box$

\noindent
{\it Proof of Lemma \ref{LemmaA4}.} Let us denote by $k(n)$ a sequence such that
 $k(n)\rightarrow 0$ and $|J_{n}v -v|_{s'}\leq k(n)|v|_{s'+1}, s'=s-2,s-1$
for all $v$ is in $H^{s'+1}$. Then we have:  
$|u_{n}-u|_{s-2}\leq |u_{n}-J_{n}u|_{s-2}+k(n)|u|_{s-1}$. Hence, to prove the lemma, it suffices to prove that  
there is $N\geq 1$  such that $ \lim_{\delta\rightarrow 0}P(\delta,n,A)=0$ where:
\[
 P(\delta,N,A):=
   {\Pr}(\sup_{t\in I}|J_{N}u(t)-u_{N}(t)|_{s-2}>
    \eta, \sup_{n}\sup_{t}|u_{n}(t)|_{s}\leq A, \sup_{t\in I}|u(t)|_{s}\leq A
    | |w|_{\infty} < \delta ).\]
\noindent
Let us set $v_{n}(t)= u_{n}(t)-J_{n}u(t)$. Then $d v_{n}(t)= a(t)v_{n}(t)\circ dw(t)+[a,J_{n}]u(t)\circ dw(t)$ and by
the It\^o formula, we have:
\begin{align*}
\PrL v_{n}(t),v_{n}(t) \PrR_{s-2}= & |J_{n}u_{0}-u_{0}|_{s-2}^{2}+ (\PrL A v_{n}(t),v_{n}(t) \PrR_{s-2}
                                  + 2 Re \PrL v_{n}(t), [a,J_{n}]u(t)\PrR_{s-2}) w(t)  \\
                                 & -\int_{0}^{t}(\PrL A v_{n}(\t),v_{n}(\t) \PrR_{s-2}+ 2 Re \PrL v_{n}(\t), [a,J_{n}]u(\t)\PrR_{s-2} )w(\t)\circ dw(\t)\\
                                 = & |J_{n}u_{0}-u_{0}|_{s-2}^{2}+  [\PrL A v_{n}(t),v_{n}(t) \PrR_{s-2}
                                 + 2 Re \PrL v_{n}(t), [a,J_{n}]u(t)\PrR_{s-2}] w(t)  \\
                                 & -\int_{0}^{t} \PrL A v_{n}(\t),v_{n}(\t) \PrR_{s-2}(\t) w(\t)\circ dw(\t)-\int_{0}^{t} X_{1}(\t) w(\t) dw(\t)\\ 
																  &- \frac{1}{2}\int_{0}^{t} X_{2}(\t) d\t,													
\end{align*}
where:
\begin{align*}
X_{1}(\t) =&  2 Re \PrL v_{n}(\t), [a,J_{n}]u(\t)\PrR_{s-2} \\
X_{2}(\t)=& w(t) [2  Re \PrL A_{n} v_{n}(\t), [a,J_{n}]u(\t)\PrR_{s-2} + 2 Re \PrL A_{n}' v_{n}(\t),  v_{n}(\t)\PrR_{s-2}]\\
					& 2 Re \PrL v_{n}(\t), [a,J_{n}]u(\t)\PrR_{s-2},
\end{align*}
\noindent
and $A_{n}'=[a,J_{n}] a  + a^{*}[a,J_{n}]$. The idea is that the terms with a multiplicative factor $w(t)$ will be controlled by the condition $|w|_{\infty}<\delta$ and
a factor $k(n)$ which tends to $0$, and the terms involving stochastic integrals will be controlled via the exponential inequality of martingales and the estimate (\ref{iw1}), but to do so we will be led to use two iterations of the It\^o formula that mimic the integration by part in stochastic calculus. So we write:
\[ \PrL v_{n}(t),v_{n}(t) \PrR_{s-2} = I_{1, n}(t)+I_{2, n}(t)+I_{3, n}(t) \]
with
\begin{align*}
I_{1, n}(t)  =& |J_{n}u_{0}-u_{0}|_{s-2}^{2}+ (\PrL A v_{n}(t),v_{n}(t) \PrR_{s-2}+ 2 Re \PrL v_{n}(t), [a,J_{n}]u(t)\PrR_{s-2}) w(t)  \\
							& - \frac{1}{2}\int_{0}^{t} X_{2}(\t) d\t		\\						
I_{2, n}(t) =& -\int_{0}^{t} \PrL A v_{n}(\t),v_{n}(\t) \PrR_{s-2}(\t) w(\t)\circ dw(\t) \\
I_{3, n}(t) =& -\int_{0}^{t} X_{1}(\t) w(\t) dw(\t) 
\end{align*}
\noindent
Using the boundedness of the operators $A$ and $L$ (and noting that for $v\in H^{s}$ we have:
$|[J_{n},a]v|_{s-2}\leq C_{1} k(n)|v|_{s}$  and $ |J_{n}v-v|_{s-2} \leq k(n)|v|_{s} $ with $k(n)\rightarrow 0$), it follows that 
on the set $\{\sup_{n}\sup_{t\in I}|u_{n}(t)|_{s}\leq A, \sup_{t\in I}|u(t)|_{s}\leq A \}$, there exist constants $C_{1}(A), C_{2}(A), C_{3}(A)$ such that
\[
I_{1, n}(t)  \leq C_{1}(A)|w|_{\infty} + k(n)  C_{2}(A).
\]
\noindent
This implies that for $n$ sufficiently large and $\delta$ sufficiently small ($n\geq N_{1}, \delta\leq \delta_{1}$), we will have $ C_{1}(A)|w|_{\infty} + k(n)  C_{2}(A) \leq \eta/2 $ and then
\begin{align*}
P(n,\delta,A) \leq & {\rm Pr}(\sup_{t\in I}(|I_{1, n}(t) | >\eta/2, \sup_{n}\sup_{ t\in I}|u_{n}(t)|_{s}\leq A,\\
		          & \sup_{t\in I}|u(t)|_{s}\leq A||w|_{\infty}<\delta) \\
		          &+ {\rm Pr}(\sup_{t\in I}(|I_{2, n}(t) | + |I_{3, n}(t)|) >\eta/2, \sup_{n}\sup_{ t\in I}|u_{n}(t)|_{s}\leq A,\\
		          &  \sup_{t\in I}|u(t)|_{s}\leq A||w|_{\infty}<\delta)\\
            \leq &   {\rm Pr}(\sup_{t\in I}|I_{3, n}(t) | >\eta/4, \sup_{n}\sup_{ t\in I}|u_{n}(t)|_{s}\leq A, \sup_{t\in I}|u(t)|_{s}\leq A||w|_{\infty}<\delta)\\
					    &+ {\rm Pr}(\sup_{t\in I}|I_{2, n}(t) | >\eta/4, \sup_{n}\sup_{ t\in I}|u_{n}(t)|_{s}\leq A, \sup_{t\in I}|u(t)|_{s}\leq A||w|_{\infty}<\delta)\\
   :=& P_{1}(\delta,n,A)+P_{2}(\delta,n,A).
\end{align*}
Using the fact that $|2 Re \PrL v_{n}(\t), [a,J_{n}]u(\t)\PrR_{s-2}|\leq k(n)C_{6}(A)$, it follows
that
\[ P_{1}(\delta,n,A)\leq C_{7}\exp-\frac{1}{C_{6}^{2}(A)k(n)^{2}\delta^{2}} C_{8}\exp \frac{C_{9}}{\delta^{2}}\]
where we have used (\ref{iw1}) and the exponential inequality (\ref{ExpIneqMartingale}).
Hence for large $N$, ($N\geq N_{2}\geq N_{1}$) we have $P_{1}(\delta,n,A)\leq C_{10}\exp{(-C_{11}/\delta^{2})}$, where $C_{10},
  C_{11} >0$ depend only on $N_{2},A$. This implies that for all
$n\geq N_{2}$ we have $\lim_{\delta\rightarrow 0}P_{1}(\delta,n,A)=0$. 
Now we shall show the same result for $P_{2}(\delta,n,A)$. By the It\^o formula we have:
\begin{align*}
I_{2, n}(t)  =&  \PrL A v_{n}(t),v_{n}(t) \PrR_{s-2} w(t)^{2} \\
                               &- \int_{0}^{t} \PrL A v_{n}(\t),v_{n}(t) \PrR_{s-2} w(\t) \circ dw(\t)\\
                           & - \int_{0}^{t}  \PrL L v_{n}(\t),v_{n}(\t) \PrR_{s-2}\circ w(\t)^{2} \circ dw(\t)\\
													& - \int_{0}^{t} 2 Re \PrL v_{n}(\t), A[a, J_{n}]u(\t) \PrR_{s-2} w(\t)^{2} \circ dw(\t)\\
													& - \frac{1}{2} \int_{0}^{t} 2 Re \PrL v_{n}(\t), A[a, J_{n}]u(\t) \PrR_{s-2} 2 w(\t) d\t  \\
													& - \frac{1}{2} \int_{0}^{t} 2 w(\t)^{2} [Re \PrL a v_{n}(\t), A[a, J_{n}]u(\t) \PrR_{s-2} 
													       + \PrL [a, J_{n}] u(\t), A[a, J_{n}]u(\t) \PrR_{s-2}\\
													&+ \PrL  v_{n}(\t), A[a, J_{n}] a u(\t) \PrR_{s-2}] d\t
\end{align*}
which yields
\begin{align*}
I_{2, n}(t)  =& \frac{1}{2}[ -\int_{0}^{t}  (\PrL L v_{n}(\t),v_{n}(t) \PrR_{s-2} + 2 Re \PrL v_{n}(\t), A[a, J_{n}]u(\t) \PrR_{s-2}) w(\t)^{2} dw(\t)\\
              &+ \PrL A v_{n}(t),v_{n}(t) \PrR_{s-2}w(t)^{2} \\
					&- \frac{1}{2}\int_{0}^{t}[ 2\PrL L v_{n}(\t),v_{n}(t) \PrR_{s-2} w(\t) +\PrL M v_{n}(\t),v_{n}(t) \PrR_{s-2} w(\t)^{2})d\t] \\
						& - \frac{1}{2}\int_{0}^{t}[  2 Re \PrL v_{n}(\t), A[a, J_{n}]u(\t) \PrR_{s-2} 2 w(\t) d\t \\
						& +  2  w(\t)^{2} Re (\PrL a v_{n}(\t), A[a, J_{n}]u(\t) \PrR_{s-2} 
							 + 2 Re \PrL [a, J_{n}] a u(\t), A[a, J_{n}]  u(\t) \PrR_{s-2}) \\
								& + \PrL av_{n}(\t), A[a, J_{n}] a u(\t) \PrR_{s-2}   ]d\t \\
									:=& I_{4,n}(t) + I_{5,n}(t),
\end{align*}
where $ I_{4,n}(t)$ is the stochastic integral term (we note in this term a factor $w(\t)^{2}$ which will be used below) and $ I_{5,n}(t)$ contains all the other terms. 
\\
\noindent
We remark that on the set $\{\sup_{n}\sup_{t\in I}|u^{n}(t)|_{s}\leq A,\sup_{t\in I}|u(t)|_{s}\leq A\}$, we have
\[ \sup_{t\in I}|I_{5,n}(t))|\leq C_{12}(A)(|w|+|w|^{2}), \] 
which implies that for $\delta$ sufficiently small, we have:
\[ P_{2}(\delta,n,A)={\rm Pr}(\sup_{t\in I}|I_{4,n}(t)|>\eta,
 \sup_{n}\sup_{t\in I}|u^{n}(t)|_{s}\leq A, \sup_{t\in I}|u(t)|_{s}\leq A  ||w|_{\infty}<\delta)   \]
\noindent
Then, since $|\PrL L v_{n}(\t),v_{n}(t) \PrR_{s-2} + 2 Re \PrL v_{n}(\t), A[a, J_{n}]u(\t) \PrR_{s-2}|$ is bounded on the set 
$\{\sup_{n}\sup_{t\in I}|u_{n}(t)|_{s}\leq A, \sup_{t\in I}|u(t)|_{s}\leq A \}$, it follows from (\ref{iw1}) and the exponential inequality that:
\[  P_{2}(\delta,n,A)\leq C_{13}\exp (-\frac{C_{14}}{\delta^{4}})
     \exp (\frac{C_{15}}{\delta^{2}}) \]
which implies that for all $n$, $\lim_{\delta\rightarrow 0} P_{2}(\delta,n,A)=
0$. This completes the proof of Lemma \ref{LemmaA4}.

\noindent
{\it Proof of Lemma \ref{LemmaA5}.} Since the operators $a_{t}(x,D)J_{n}$ and $a_{t}(x,D)J_{n}$ are bounded, the proof of this lemma can be done exactly as in
the finite-dimensional case, see Ikeda-Watanabe \cite{IW}, Theorem 8.2, p.419.

%%%%%%%%%%%%%%%%%%%%%%%%%%%%%%%%%%%%%%%%%%%%%%%%%%%%%%%%%%%%%%%%%%%%%%5

\textsc{Acknowledgments.} I am grateful to the anonymous Referee for valuable remarks, corrections and suggestions that lead to clarify the paper and for helpful advice about its presentation.
\\
{\small 

BP. 2119 Rabat Hay Riyad, \\
Rabat, \\
Morocco\\
e-mail: adnan.aboulalaa@polytechnique.org
\end{document}